\documentclass[runningheads]{llncs}
\usepackage{graphicx}
\usepackage{booktabs}
\usepackage{multirow}

\usepackage{hyperref}

\usepackage{enumitem}

\usepackage{subfig}
\usepackage{mwe}

\usepackage{pifont}

\usepackage{mathtools} 
\usepackage{amsfonts}
\usepackage{bm} 

\newcommand{\argmin}{\operatornamewithlimits{argmin}}

\usepackage{amsmath} 
\usepackage{algorithm} 
\usepackage{algpseudocode}
\makeatletter
\renewcommand{\ALG@beginalgorithmic}{\scriptsize}
\makeatother

\usepackage{nccmath} 
\usepackage{etoolbox}
\newenvironment{myequation}{%
\smallskip\par\centering$\displaystyle}
{$\smallskip\par}

\usepackage[square,sort,comma,numbers]{natbib} 

\begin{document}
\title{Multi-Task Predict-then-Optimize}
%
%
\author{Bo Tang\orcidID{0000-0002-6035-5167} \and
Elias B.~Khalil\orcidID{0000-0001-5844-9642}}

\authorrunning{B. Tang \and E. Khalil}
%
\institute{SCALE AI Research Chair in Data-Driven Algorithms for Modern Supply Chains\\Department of Mechanical and Industrial Engineering, University of Toronto\\
\email{\{botang,khalil\}@mie.utoronto.ca}}
\maketitle              
\begin{abstract}
The predict-then-optimize framework arises in a wide variety of applications where the unknown cost coefficients of an optimization problem are first predicted based on contextual features and then used to solve the problem. In this work, we extend the predict-then-optimize framework to a multi-task setting: contextual features must be used to predict cost coefficients of multiple optimization problems, possibly with different feasible regions, simultaneously. For instance, in a vehicle dispatch/routing application, features such as time-of-day, traffic, and weather must be used to predict travel times on the edges of a road network for multiple traveling salesperson problems that span different target locations and multiple $s-t$ shortest path problems with different source-target pairs. We propose a set of methods for this setting, with the most sophisticated one drawing on advances in multi-task deep learning that enable information sharing between tasks for improved learning, particularly in the small-data regime. Our experiments demonstrate that multi-task predict-then-optimize methods provide good tradeoffs in performance among different tasks, particularly with less training data and more tasks. 


\end{abstract}
\section{Introduction} \label{sec:intro}

The predict-then-optimize framework, in which the unknown coefficients for an optimization problem are predicted and then used to solve the problem, is emerging as a useful framework in some applications. For instance, in vehicle routing and job scheduling, we often require optimization where the model's cost coefficients, e.g., travel time and execution time, are unknown but predictable at decision time. In the conventional two-stage method, a learning model is first trained to predict cost coefficients, after which a solver separately optimizes accordingly. However, end-to-end approaches that learn predictive models that minimize the decision error directly have recently gained interest due to some improvements in experimental performance~\cite{bengio1997using, ford2015beware, donti2017task, elmachtoub2021smart, mandi2020interior, mandi2020smart}. Although there has been some recent work in predicting elements of the constraint matrix in a linear programming setting~\cite{estes2023smart,hu2023predict}, our focus here is on the predominant line of research that has focused on unknown cost coefficients~\cite{donti2017task, wilder2019melding, ferber2020mipaal, mandi2020interior, agrawal2019differentiable, elmachtoub2021smart, mandi2020smart, poganvcic2019differentiation, berthet2020learning, dalle2022learning}.

Previous work on predict-then-optimize has focused on learning the cost coefficients for a single optimization task. However, it is natural to consider the setting where multiple related tasks can share information and representations. For example, a vehicle routing application requires predicting travel times on the edges of a road network for multiple traveling salesperson problems (TSPs) that span different target locations and multiple shortest path problems with different source-target pairs. These travel time predictions should be based on the~\textit{same contextual information}, e.g., if the tasks are to be executed at the same time-of-day, then the travel times that should be predicted for the different tasks depend on the same features. Another case is in package delivery, where distributing packages from one depot to multiple depots results in independent delivery tasks that nonetheless share travel time predictions since they  use the same road network. To that end, we introduce multi-task end-to-end predict-then-optimize, which simultaneously solves multiple optimization problems with a loss function that relates to the decision errors of all such problems.

Multi-task learning has been successfully applied to natural language processing, computer vision, and recommendation systems. However, its applicability to the predict-then-optimize paradigm is yet to be explored. Predict-then-optimize with multi-task learning is attractive because of the ability to improve model performance in the small-data regime. Machine learning, especially with deep neural networks, is data-intensive and prone to overfitting, which might limit applicability to the predict-then-optimize paradigm. The need to simultaneously minimize the losses of different tasks helps reduce overfitting and improve generalization. Multi-task learning combines the data of all tasks, which increases the overall training data size and alleviates task-specific noise. 

To the best of our knowledge, we introduce multi-task learning for end-to-end predict-then-optimize for the first time. We motivate and formalize this problem before proposing a set of methods for combining the different task losses. Our experiments show that multi-task end-to-end predict-then-optimize provides performance benefits and good tradeoffs in performance among different tasks, especially with less training data and more tasks, as compared to applying standard single-task methods independently for each task. As an additional contribution, we distinguish end-to-end predict-then-optimize approaches that learn from observed costs (the usual setting of~\citet{elmachtoub2021smart}) and those that learn directly from (optimal) solutions without the objective function costs themselves. This extends our framework to applications where there are no labeled coefficients in the training data, e.g., the Amazon Last Mile Routing Challenge~\cite{winkenbach2021technical}. The open-source code is available \footnote{\url{https://github.com/khalil-research/Multi-Task_Predict-then-Optimize}}.
\section{Related Work} \label{sec:liter}

\subsection{Differentiable Optimization}
\label{subsec:diff}

The key component of gradient-based end-to-end predict-then-optimize is differentiable optimization, which allows the backpropagation algorithm to update model parameters from the decision made by the optimizer. Based on the KKT conditions,~\citet{amos2017optnet} introduced OptNet, a differentiable optimizer to make quadratic optimization problems learnable. With OptNet,~\citet{donti2017task} investigated a learning framework for quadratic programming;~\citet{wilder2019melding} then added a small quadratic regularization to the linear objective function for linear programming;~\citet{ferber2020mipaal} extended the method to the discrete model with the cutting plane;~\citet{mandi2020interior} adopted log-barrier regularization instead of a quadratic one. Besides OptNet, \citet{agrawal2019differentiable} leveraged KKT conditions to differentiate conic programming.

Except for the above approaches with KKT, an alternative methodology is to design gradient approximations to avoid ill-defined gradients from predicted costs to optimal solutions.~\citet{elmachtoub2021smart} proposed a convex surrogate loss.~\citet{poganvcic2019differentiation} developed a differentiable optimizer through implicit interpolation.~\citet{berthet2020learning} demonstrated a method with stochastic perturbation to smoothen the loss function and further constructed the Fenchel-Young loss.~\citet{dalle2022learning} extended the perturbation approach to the multiplicative perturbation and the Frank-Wolfe regularization.~\citet{mulamba2020contrastive} studied a solver-free contrastive loss. Moreover,~\citet{shahdecision} provided an alternate paradigm that additionally trains a model to predict decision errors to replace the solver.

\subsection{Multi-Task Learning}
\label{mtl}

Multi-task learning, first proposed by~\citet{caruana1997multitask}, aims to learn multiple tasks together with joint losses. In short, a model with multiple loss functions is defined as multi-task learning. Much research has focused on the (neural network) model architecture: the most basic model is a shared-bottom model~\citep{caruana1997multitask}, including shared hidden layers at the bottom and task-specific layers at the top. Besides such hard parameter sharing schemes, there is also soft sharing so that each task keeps its own parameters.~\citet{duong2015low} added $l_2$ norm regularization to encourage similar parameters between tasks. Furthermore, neural networks with different experts and some gates~\cite{shazeer2017outrageously, ma2018modeling, tang2020progressive} were designed to fuse information adaptively.

Another crucial issue is resolving the unbalanced (in magnitude) and conflicting gradients from different tasks. There are weighting approaches, such as UW~\cite{kendall2018multi}, GradNorm~\cite{chen2018gradnorm}, and DWA~\cite{liu2019end}, that have been proposed to adjust the weighting of different losses. Other methods, such as PCGrad~\cite{yu2020gradient}, GradVec~\cite{wang2020gradient}, and CAGrad~\cite{liu2021conflict}, were designed to alter the direction of the overall gradient to avoid conflicts and accelerate convergence.
\section{Building Blocks} \label{sec:block}

\subsection{Optimization Problem}
\label{subsec:def}

For multi-task end-to-end predict-then-optimize, each task $t$ is a separate (integer) linear optimization problem, defined as follows:

\begin{myequation}
\begin{aligned}
\min_{\bm{w}^t} \quad & {\bm{c}^t}^T \bm{w}^t \\
\textrm{s.t.} \quad & \bm{A}^t \bm{w}^t \leq \bm{b}^t \\
& \bm{w}^t \geq \bm{0} \\
& \textrm{Some } w^t_i \textrm{ are integer.}
\end{aligned}
\end{myequation}

The decision variables are $\bm{w}^t$, the constraint coefficients are $\bm{A}^t$, the right-hand sides of the constraints are $\bm{b}^t$, and the unknown cost coefficients are $\bm{c}^t$. Given a cost $\bm{c}^t$, $\bm{w}^{t*}_{\bm{c}^t}$ is the corresponding optimal solution, and $z^{t*}_{\bm{c}^t}$ is the optimal  value. Additionally, we define $S^t = \{\bm{A}^t \bm{w}^t \leq \bm{b}^t, \bm{w}^t \geq \bm{0}, ...\}$ as the feasible region of decision variables $\bm{w}^t$.

\subsection{Gradient-based Learning}
\label{subsec:e2e}

\begin{figure}[htbp]
    \centering
    \includegraphics[width=0.9\textwidth]{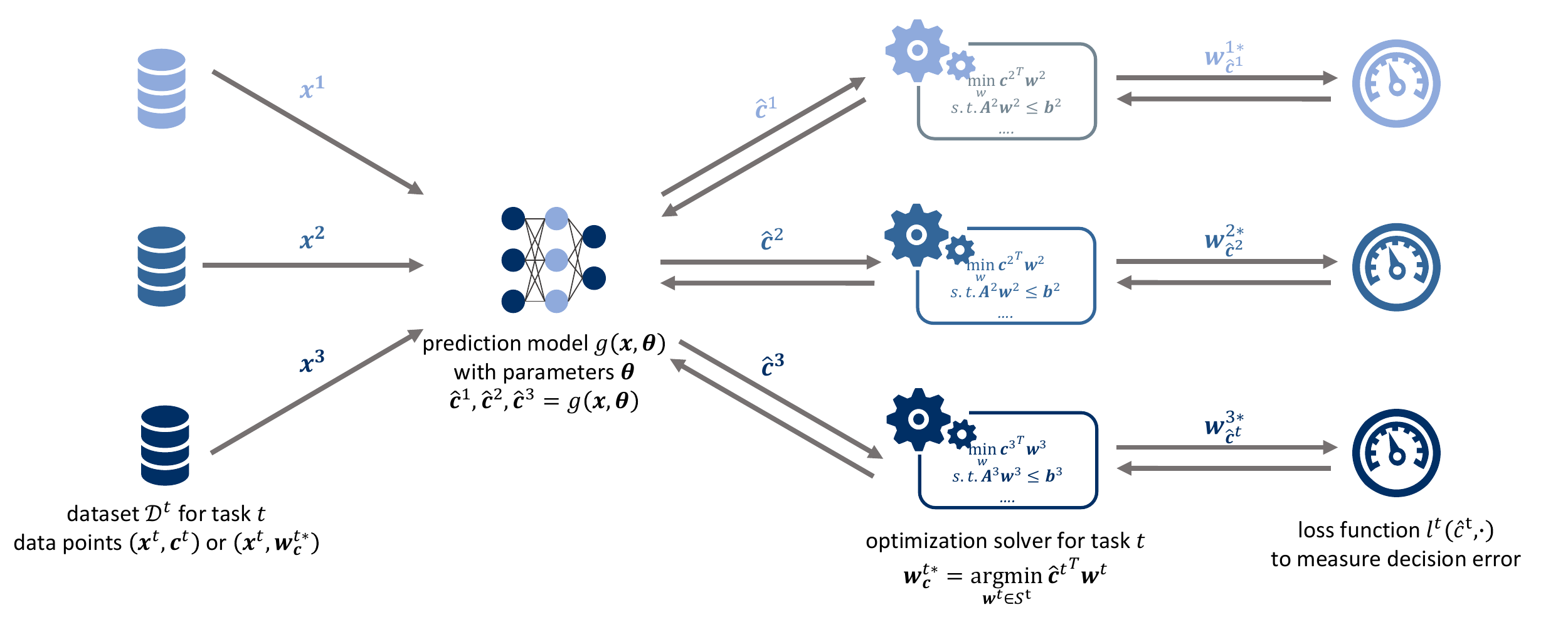}    
    \caption{Illustration of end-to-end multi-task end-to-end predict-then-optimize: labeled datasets $\mathcal{D}^1, \mathcal{D}^2, ..., \mathcal{D}^t$ are used to fit a machine learning predictor $g(\bm{x}; \bm{\theta})$ that predicts costs $\hat{\bm{c}}^t$ for each task $t$. The loss function $l^t (\hat{\bm{c}}^t, \cdot)$ to be minimized measures decision error instead of prediction error.}
    \label{fig:multi}
\end{figure}

End-to-end predict-then-optimize aims to minimize a decision loss directly. As a supervised learning problem, it requires a labeled dataset $\mathcal{D}$ consisting of features $\bm{x}$ and labels $\bm{c}$ or $\bm{w}^*_c$. We will further discuss the difference between cost labels $\bm{c}$ and solution labels $\bm{w}^*_c$ in Sec \ref{subsec:label}. As shown in Figure \ref{fig:multi}, multi-task end-to-end predict-then-optimize predicts the unknown costs for multiple optimization problems and then solves these optimization problems with the predicted costs. The critical component is a differentiable optimization embedded into a differentiable predictive model. However, the only learnable part is the prediction model $g(\bm{x}; \bm{\theta})$, since there are no parameters to update in the solver and loss function.


\subsection{Decision Losses}
\label{subsec:loss}

Loss functions for end-to-end predict-then-optimize aim to measure the error in decision-making. For instance, regret is defined as the difference in objective values due to an optimal solution that is based on the true costs and another that is based on the predicted costs:
\begin{myequation}
    l_{\text{Regret}}(\hat{\bm{c}}, \bm{c}) = \bm{c}^T \bm{w}^*_{\hat{\bm{c}}} - z^*_{\bm{c}}
    \label{eq:regret}
\end{myequation}

However, with a linear objective function, regret does not provide useful gradients for learning~\cite{elmachtoub2021smart}. Besides regret, decision error can also be defined as the difference between the true solution and its prediction, such as using the Hamming distance~\cite{poganvcic2019differentiation} or the squared error of a solution to an optimal one~\cite{berthet2020learning}. 
Because the function from costs $\bm{c}$ to optimal solutions $\bm{w}^*_{\bm{c}}$ is piecewise constant, a solver with any of the aforementioned losses has no nonzero gradients to update the model parameter $\bm{\theta}$. Thus, the state-of-art methods, namely Smart Predict-then-Optimize (SPO+)~\cite{elmachtoub2021smart} and Perturbed Fenchel-Young Loss (PFYL), both design surrogate decision losses which allow for a nonzero approximate gradient (or subgradient), $\frac{\partial l(\cdot)}{\partial \hat{\bm{c}}}$.\\

\noindent\textbf{Smart Predict-then-Optimize Loss (SPO+).} 
SPO+ loss~\cite{elmachtoub2021smart} is a differentiable convex upper bound on the regret: 
\begin{myequation}
\begin{aligned}
l_{\text{SPO+}}(\hat{\bm{c}}, \bm{c}) = - \underset{\bm{w} \in S}{\min} \{(2 \hat{\bm{c}} - \bm{c})^T \bm{w}\} + 2 \hat{\bm{c}}^T \bm{w}^*_{\bm{c}} - z^*_{\bm{c}}.
\end{aligned}
\end{myequation} 
\noindent One proposed subgradient for this loss writes as follows:
\begin{myequation}
\begin{aligned}
2 (\bm{w}^*_{\bm{c}} - \bm{w}^*_{2 \hat{\bm{c}} - \bm{c}}) \in \frac{\partial l_{\text{SPO+}}(\hat{\bm{c}}, \bm{c})}{\partial \hat{\bm{c}}}.
\end{aligned}
\end{myequation}

In theory, the SPO+ framework can be applied to any problem with a linear objective.~\citet{elmachtoub2021smart} have conducted experiments on shortest path as a representative linear program and portfolio optimization as a representative quadratically-constrained problem.\\

\noindent\textbf{Perturbed Fenchel-Young Loss (PFYL).} 
  PFYL~\cite{berthet2020learning} leverages Fenchel duality, where \citet{berthet2020learning} discussed only the case of linear programming. The predicted costs are sampled with Gaussian perturbation $\bm{\xi}$, and the expected function of the perturbed minimizer is defined as $F(\bm{c}) =  \mathop{{}\mathbb{E}}_{\bm{\xi}}[\underset{\bm{w} \in S}{\min} {\{(\bm{c}+\sigma \bm{\xi})^T \bm{w}\}}]$. With $\Omega (\bm{w}^*_{\bm{c}}),$ the dual of $F(\bm{c})$, the Fenchel-Young loss reads:
$$l_{\text{FY}}(\hat{\bm{c}}, \bm{w}^*_{\bm{c}}) =  \hat{\bm{c}}^T \bm{w}^*_{\bm{c}} - F(\hat{\bm{c}}) - \Omega (\bm{w}^*_{\bm{c}}).
$$
Then, we can estimate the gradients by $M$ samples Monte Carlo: $$\frac{\partial l_{\text{FY}}(\hat{\bm{c}}, \bm{w}^*_{\bm{c}})}{\partial \hat{\bm{c}}} \approx \bm{w}^*_{\bm{c}} - \frac{1}{M} \sum_m^M \underset{\bm{w} \in S}{\argmin} {\{(\hat{\bm{c}}+\sigma \bm{\xi}_m)^T \bm{w}\}}$$

\subsection{Multi-Task Loss Weighting Strategies}
\label{subsec:multi}

The general idea of multi-task learning is that multiple tasks are solved simultaneously by the same predictive model. It is critical for a multi-task neural network, one such flexible class of models, to balance losses among tasks with the loss weights $u_t$ for task $t$. The weighting approaches we evaluated include a uniform combination (all $u^t = 1$) and GradNorm~\cite{chen2018gradnorm}, an adaptive loss weighting approach. The latter provides adaptive per-task weights $u_{\text{Ada}}^t$ that are dynamically adjusted during training in order to keep the scale of the gradients similar. In this work, we set the GradNorm hyperparameters of ``restoring force" to $0.1$ and the learning rate of loss weights to $0.005$. Further tuning is possible but was not needed for our experiments.

All the training strategies we have explored in this paper, including baseline approaches, are summarized in Table~\ref{tab:loss}. 
Let $T$ be the number of tasks, and cost coefficient prediction for task $t$ be $\hat{\bm{c}}^t$.
``mse" is the usual two-stage baseline of training a regression model that minimizes cost coefficients mean-squared error $l_{\text{MSE}} = \frac{1}{n} \sum_i^n \| \hat{\bm{c}}_i - \bm{c}_i \|^2$ only without regard to the decision. ``separated" trains one model per task, minimizing, for each task, a decision-based loss such as SPO+ or PFYL from Section~\ref{subsec:loss}. ``comb" simply sums up the per-task decision losses, whereas ``gradnorm" does so in a weighted adaptive way. For any of these methods, whenever ``+mse" is appended to the method name, a variant of the method is obtained that combines additional mean-squared error $l_{\text{MSE}}$ in the cost predictions with the decision loss. Such a regularizer is known to be useful in practice, even with the primary evaluation metric of a trained model being its decision regret~\cite{elmachtoub2021smart}.
Although we refer to ``separated+mse'' as a single-task method, it can also be considered as a multi-task learning method in a broad sense because of the inclusion of two losses. 


\begin{table}[t]
\centering
\begin{tabular}{l|l|l}
\hline
& \textbf{Strategy} & \textbf{Losses}                                                                                                                        \\ \hline
\multirow{3}{*}{Single-Task}
& mse               & $l_{\text{MSE}} (\bm{c}, \hat{\bm{c}})$                                                                                                \\
& separated         & Separate $l_{\text{Decision}} (\hat{\bm{c}}^t, \cdot)$ for each task $t$                                                                               \\
& separated+mse     & Separate $l_{\text{Decision}} (\hat{\bm{c}}^t, \cdot) + l_{\text{MSE}} (\bm{c}^t, \hat{\bm{c}}^t)$  for each task $t$                                      \\ \hline 
\multirow{4}{*}{Multi-Task}
& comb              & $\sum_t^T l_{\text{Decision}}^t (\hat{\bm{c}}^t, \cdot)$                                                                                               \\
& comb+mse          & $\sum_t^T l_{\text{Decision}}^t (\hat{\bm{c}}^t, \cdot) + \sum_t^T l_{\text{MSE}}^t (\bm{c}^t, \hat{\bm{c}}^t)$                                            \\
& gradnorm          & $\sum_t^T u_{\text{Ada}}^t l_{\text{Decision}}^t (\hat{\bm{c}}^t, \cdot)$                                                                              \\
& gradnorm+mse      & $\sum_t^T u_{\text{Ada}_1}^t l_{\text{Decision}}^t (\hat{\bm{c}}^t, \cdot) + \sum_t^T u_{\text{Ada}_2}^t l_{\text{MSE}}^t (\bm{c}^t, \hat{\bm{c}}^t)$          \\  \hline        
\end{tabular}
\caption{Losses of Various Training Strategies}\label{tab:loss}
\end{table}

\section{Learning Architectures} \label{sec:mthd}

\subsection{Shared Learnable Layers}
\label{subsec:share}

\begin{figure}[htbp]
    \centering
    \includegraphics[width=0.75\textwidth]{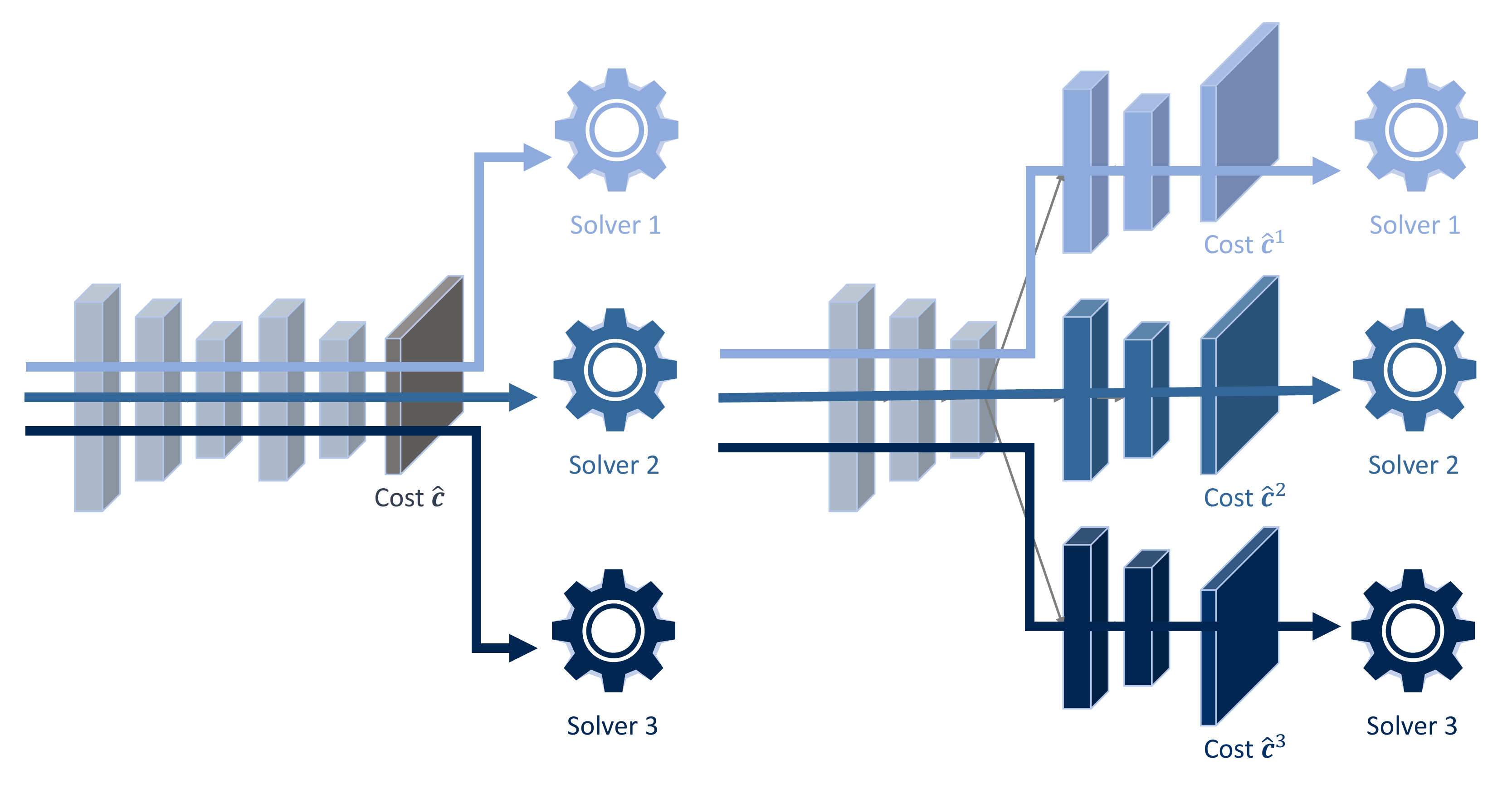}    
    \caption{Illustration of two types of multi-task end-to-end predict-then-optimize: On the left, all optimization tasks require the same prediction as cost coefficients. On the right, different tasks share some layers as feature embeddings and make different coefficient predictions.}
    \label{fig:cost}
\end{figure}

The model class we will explore is deep neural networks. Besides their capacity to represent complex functions from labeled data, neural networks have a compositional structure that makes them particularly well-suited for multi-task learning. A multi-task neural network shares hidden layers across all tasks and keeps specific layers for each task. Figure \ref{fig:cost} illustrates that the sharing part of multi-task end-to-end predict-then-optimize depends on the consistency of the predicted coefficients, which we will define next. At a high level, Figure~\ref{fig:cost} distinguishes two settings. On the left, the different tasks use the exact same predicted cost vector. On the right, each task could have a different cost vector. In both settings, the predictions are based on the same input feature vector.\\ 

\noindent\textbf{Shared Predicted Coefficients (Single-Cost).} 
In this setting (left of Figure \ref{fig:cost}), which we will refer to as~\textit{single-cost}, the same cost coefficients are shared among all tasks. For example, multiple navigation tasks on a single map are shortest-path problems with different source-target pairs that share the same distance matrix (i.e., costs). In this case, the cost coefficients $\bm{c}^t$ for task $t$ are equal to or a subset of the shared costs $\bm{c}$. Thus, the prediction model is defined as $$\hat{\bm{c}} = g(\bm{x}; \bm{\theta}_{\text{Shared}}),$$ which is the same as a single-task model, and the multiple tasks combine their losses, $$\sum_t u^t l_{\text{Decision}}^t (\hat{\bm{c}}, \cdot),$$ based on the shared prediction, $\hat{\bm{c}}$. Therefore, as Algorithm~\ref{alg:single} shows, all learnable layers are shared. In addition, the baseline methods we referred to as ``separated'' and ``separated+mse'' are not  practical in this same-costs setting as they, inconsistently, produce different cost predictions for each task, even when that is not required. Nonetheless, an experimental assessment of their performance will be conducted in order to contrast it with using multi-task learning.\\

\begin{algorithm}
  \caption{Single-Cost Multi-Task Gradient Descent}\label{alg:single}
  \small
  \begin{algorithmic}[1]
    \Require coefficient matrix $\bm{A}^1, \bm{A}^2, ...$; right-hand side $\bm{b}^1, \bm{b}^2, ...$; training data $\mathcal{D}$
    \State Initialize predictor parameters $\bm{\theta}_{\text{Shared}}$ for predictor $g(\bm{x}; \bm{\theta}_{\text{Shared}})$
    \For{epochs}
      \For{each batch of training data $(\bm{x}, \bm{c})$ or $(\bm{x}, \bm{w}_{\bm{c}}^{1*}, \bm{w}_{\bm{c}}^{2*}, ...)$}
        \State Forward pass to predict cost $\hat{\bm{c}} := g(\bm{x}; \bm{\theta}_{\text{Shared}})$
        \State Forward pass to solve optimal solution $\bm{w}^{t*}_{\hat{\bm{c}}} :=  {\argmin}_{\bm{w}^t \in S^t} \hat{\bm{c}}^T \bm{w}^t$ per task
        \State Forward pass to sum weighted decision losses $l(\hat{\bm{c}}, \cdot) := \sum_t u^t l^t (\hat{\bm{c}}, \cdot)$ 
        \State Backward pass from loss $l(\hat{\bm{c}}, \cdot)$  to update parameters $\bm{\theta}_{\text{Shared}}$
      \EndFor
    \EndFor
  \end{algorithmic}
\end{algorithm}

\noindent\textbf{Shared Features Embeddings (Multi-Cost).} 
In many applications of predict-then-optimize, the optimization problem requires cost coefficients that are specific and heterogeneous to each task, but that can be inferred from homogeneous contextual features. For instance, in a vehicle routing application, features such as time of day and weather predict travel time in different regions. Compared to the single-cost setting we just introduced, the~\textit{multi-cost} predictor here has the form $$\hat{\bm{c}}^t = g(\bm{x}; \bm{\theta}_{\text{Shared}}; \bm{\theta}_t).$$ Per-task predictions are made by leveraging the same information embedding in the layers of the neural network that are shared across the tasks (see right of Figure~\ref{fig:cost}). Thus, the corresponding loss function is $$\sum_t u^t l_{\text{Decision}}^t (\hat{\bm{c}}^t, \cdot),$$ and Algorithm~\ref{alg:multi} updates the parameters of the predictor.

\begin{algorithm}
  \caption{Multi-Cost Multi-Task Gradient Descent}\label{alg:multi}
  \small
  \begin{algorithmic}[1]
    \Require coefficient matrix $\bm{A}^1, \bm{A}^2, ...$; right-hand side $\bm{b}^1, \bm{b}^2, ...$; training data $\mathcal{D}^1, \mathcal{D}^2, ...$
    \State Initialize predictor parameters $\bm{\theta}_{\text{Shared}}, \bm{\theta}_1, \bm{\theta}_2, ...$ for predictor $g(\bm{x}; \bm{\theta}_{\text{Shared}}; \bm{\theta}_t)$
    \For{epochs}
      \For{each batch of training data $(\bm{x}^1, \bm{x}^2, ..., \bm{c}^1, \bm{c}^2, ...)$ or $(\bm{x}^1, \bm{x}^2, ..., \bm{w}_{\bm{c}^1}^{1*}, \bm{w}_{\bm{c}^2}^{2*}, ...)$}
        \State Forward pass to predict cost $\hat{\bm{c}}^t := g(\bm{x}^t; \bm{\theta}_{\text{Shared}}; \bm{\theta}_t)$ per task
        \State Forward pass to solve optimal solution $\bm{w}^{t*}_{\hat{\bm{c}}^t} :=  {\argmin}_{\bm{w}^t \in S^t} \hat{\bm{c}}^{tT} \bm{w}^t$ per task
        \State Forward pass to sum weighted decision losses $l(\hat{\bm{c}}^1, \hat{\bm{c}}^2, ..., \cdot) := \sum_t u^t l^t (\hat{\bm{c}}^t, \cdot)$ 
        \State Backward pass from loss $l(\hat{\bm{c}}^1, \hat{\bm{c}}^2, ..., \cdot)$  to update parameters $\bm{\theta}_{\text{Shared}}, \bm{\theta}_1, \bm{\theta}_2, ...$
      \EndFor
    \EndFor
  \end{algorithmic}
\end{algorithm}

\subsection{Label Accessibility and Learning Paradigms}
\label{subsec:label}

We distinguish two learning paradigms that require different kinds of labels (cost coefficients $\bm{c}$ or optimal solutions $\bm{w}^*_{\bm{c}}$) in the training data: \textbf{learning from costs} and \textbf{learning from (optimal) solutions}. This distinction is based on the availability of labeled cost coefficients $\bm{c}$. Thus, SPO+ is learning from costs because the calculation of SPO+ loss involves true cost coefficients $\bm{c}$, whereas PFYL is learning from solutions that do not require access to $\bm{c}$; see Section~\ref{subsec:loss}.

The need for the true cost coefficients as labels in the training data is a key distinguishing factor because these cost coefficients provide additional information that can be used to train the model, but they may be absent in the data. Deriving optimal solutions from the cost coefficients is trivial, but the opposite is intricate as it requires some form of inverse optimization~\cite{ahuja2001inverse}. The ability to directly learn from solutions extends the applicability of end-to-end predict-then-optimize beyond what a two-stage approach, which is based on regressing on the cost coefficients, can do. Indeed, the recent MIT-Amazon Last Mile Routing Challenge~\cite{winkenbach2021technical} is one such example in which good TSP solutions are observed on historical package delivery instances, but the corresponding edge costs are unobserved. Those good solutions are based on experienced drivers' tacit knowledge.
\section{Experiments} \label{sec:eval}

In this section, we present experimental results for multi-task end-to-end predict-then-optimize. In our experiments, we evaluate decision performance using regret, and we use mean-squared error (MSE) to measure the prediction of cost coefficients $\hat{\bm{c}}$. We use SPO+ and PFYL as typical methods for learning from costs and learning from solutions and adopt various multi-task learning strategies discussed in Sec \ref{subsec:multi}, as well as two-stage and single-task baselines. Our experiments are conducted on two datasets, including graph routing on  PyEPO TSP dataset~\footnote{\url{https://khalil-research.github.io/PyEPO}}~\cite{paszke2019pytorch}, and adjusted Warcraft terrain~\footnote{ \url{https://drive.google.com/file/d/1lYPT7dEHtH0LaIFjigLOOkxUq4yi94fy}}~\cite{poganvcic2019differentiation} to learn single-cost decisions and multi-cost decisions. We also vary the amount of training data size and the number of tasks.

All the numerical experiments were conducted in Python v3.7.9 with two Intel E5-2683 v4 Broadwell CPUs, two NVIDIA P100 Pascal GPUs, and 8GB memory. Specifically, we used PyTorch~\cite{paszke2019pytorch} v1.10.0 for the prediction model and Gurobi~\cite{gurobi} v9.1.2 for the optimization solver, and PyEPO~\cite{tang2022pyepo} v0.2.0 for SPO+ and PFYL autograd functions.

\subsection{Benchmark Datasets and Neural Network Architecture}
\label{data}

\noindent\textbf{Graph Routing with Multiple Tasks.} 
We used the traveling salesperson problem dataset generated from PyEPO~\cite{tang2022pyepo}, which uses the Euclidean distance among nodes plus polynomial function $f(\bm{x}_i) = (\frac{1}{\sqrt{p}} (\mathcal{B} \bm{x}_i)_j + 3)^4$ (where $\mathcal{B}$ is a random matrix) with random noise perturbations $f(\bm{x}_i) \cdot \epsilon$ to map the features $\bm{x}$ into a symmetric distance matrix of a complete graph. We discuss both learning from costs and learning from solutions. In this experiment, the number of features $p$ is $10$, the number of nodes $m$ is $30$, the polynomial degree of function $f(\bm{x}_i)$ is $4$, and the noise $\epsilon$ comes from U(0.5, 1.5). We sample $15 - 22$ nodes as target locations for multiple traveling salesperson problems (TSPs) and $54$ undirected edges for multiple shortest paths (SPs) with different source-target pairs. Thus, all TSP and SP tasks share the same cost coefficients.

Since the multiple routing tasks require consistent cost coefficients, the model $g(\bm{x}; \bm{\theta})$ makes one prediction of the costs that is used for all of the tasks. The architecture of the regression network is one fully-connected layer with a softplus activation to prevent negative cost predictions, and all tasks share the learnable layer. For the hyperparameters, the learning rate is $0.1$, the batch size is $32$, and the max training iterations is $30000$ with $5$ patience early stopping. For PFYL, the number of samples $M$ is $1$, and the perturbation temperature $\sigma$ is $1.0$. We formulate  SP as a network flow problem and use the Dantzig–Fulkerson–Johnson (DFJ) formulation \cite{dantzig1954solution} to solve TSP.\\

\begin{figure}[htbp]
    \centering
    \includegraphics[width=1.0\textwidth]{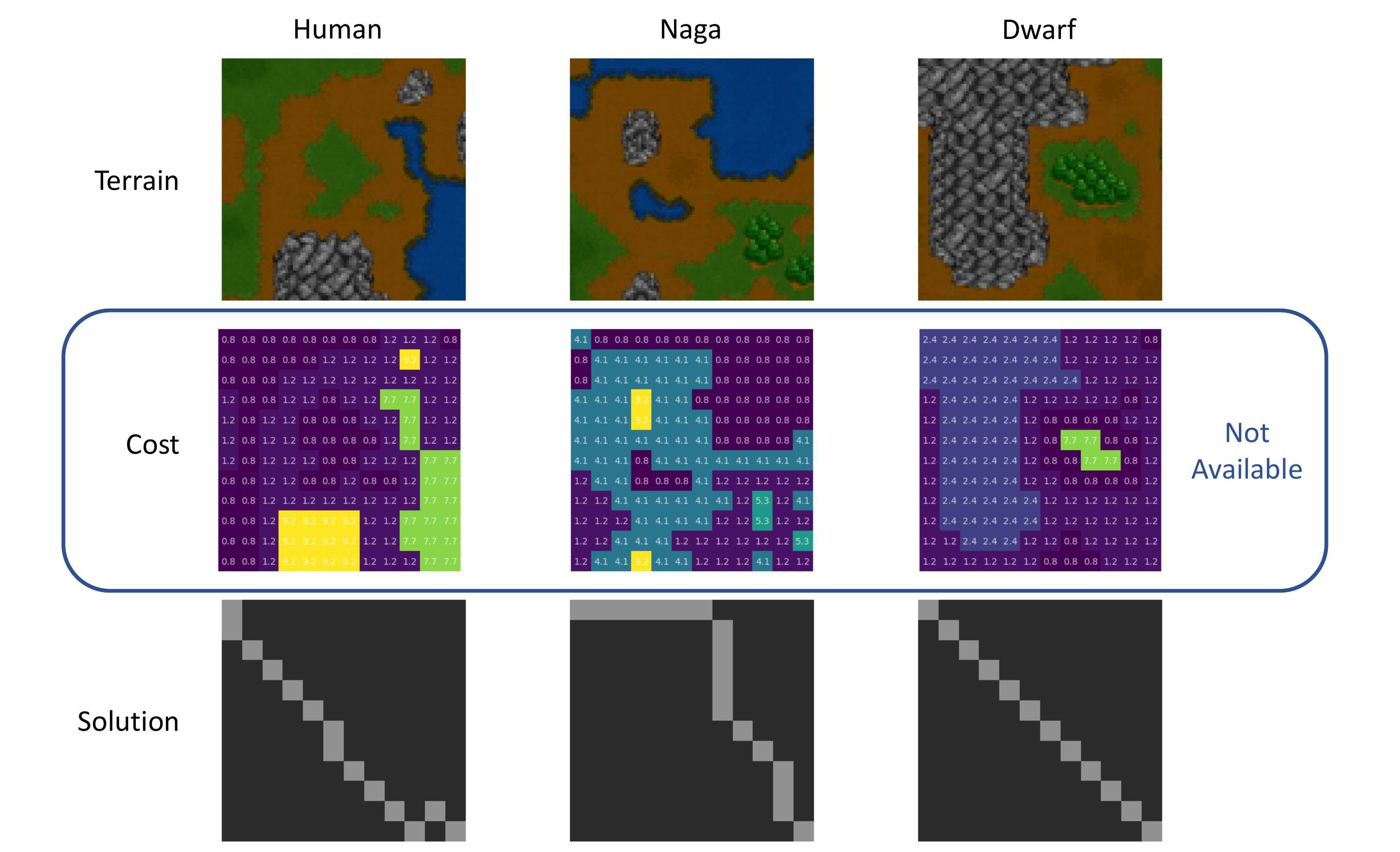}    
    \caption{Multiple datasets of Warcraft terrain images for different species, where labeled cost coefficients are unavailable.}
    \label{fig:wc}
\end{figure}

\noindent\textbf{Warcraft Shortest Path with Various Species.} 
The Warcraft map shortest path dataset~\cite{poganvcic2019differentiation} allows for the learning of the shortest path from RGB terrain images, and we use $96 \times 96$ RGB images for $12 \times 12$ grid networks and sample $3$ small subsets from $10000$ original training data points for use in training. As shown in Figure \ref{fig:wc}, we modify the cost coefficients for different species (human, naga, dwarf) and assume that the cost coefficients are not accessible in the data. This means there are three separate datasets of feature-solution pairs, which require learning from solutions using the PFYL method. Similar to SP tasks in Graph Routing, the shortest path optimization model is a linear program.

Since the multiple Warcraft shortest paths tasks require us to predict cost coefficients for different species, the prediction model should incorporate task-specific layers. Following~\citet{poganvcic2019differentiation}, we train a truncated ResNet18 (first five layers) for $50$ epochs with batches of size $70$, and learning rate $0.0005$ decaying at the epochs $30$ and $40$. The first three layers are the shared-bottom. The number of samples $M$ is $1$, and the perturbation temperature $\sigma$ is $1.0$. 

\subsection{Performance Advantage of Multi-Task Learning}
\label{subsec:perform}

\begin{figure}[htbp]
    \centering
    \subfloat{\includegraphics[width=0.50\textwidth]{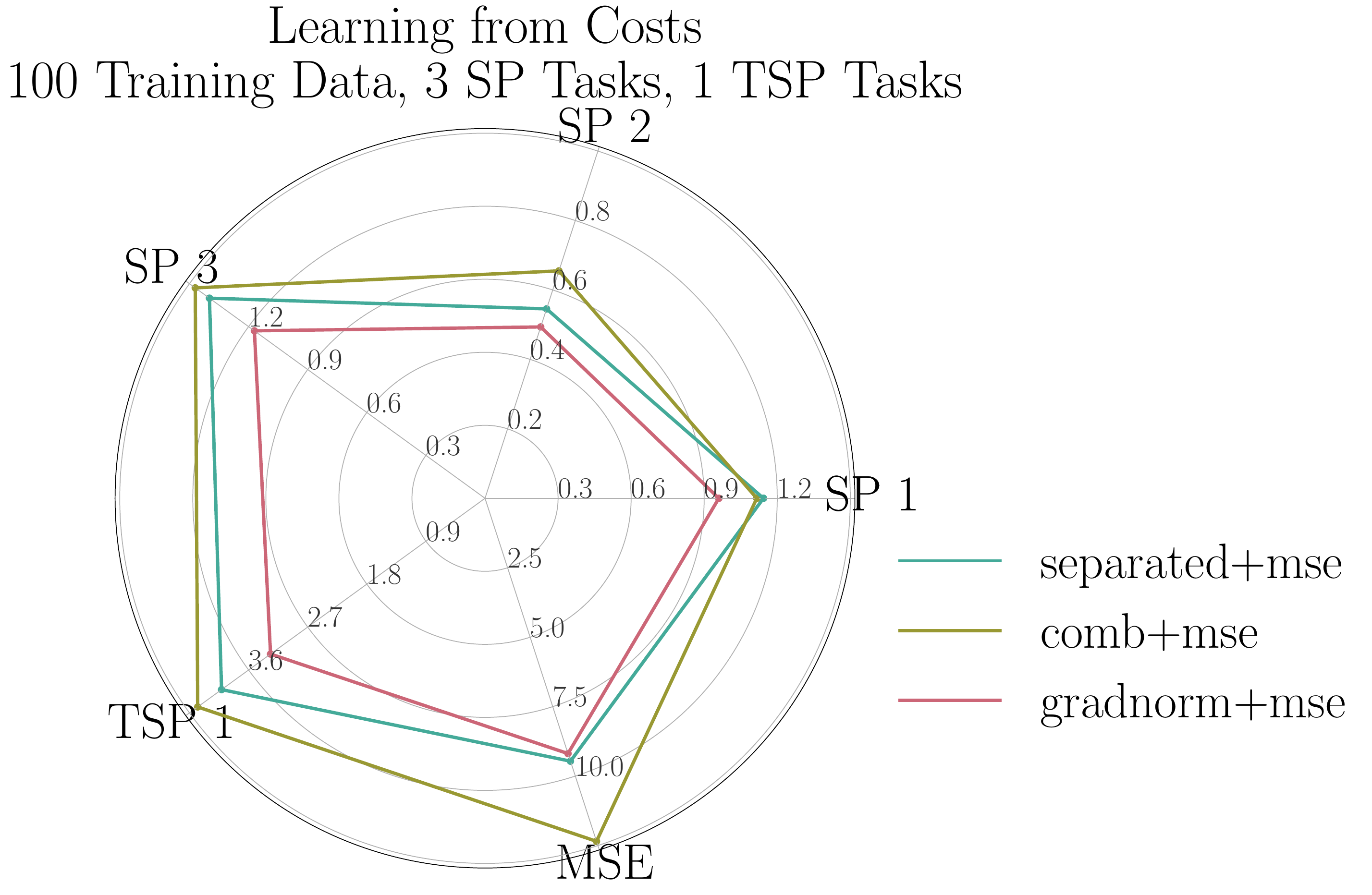}} \hspace{-0.2cm}
    \subfloat{\includegraphics[width=0.50\textwidth]{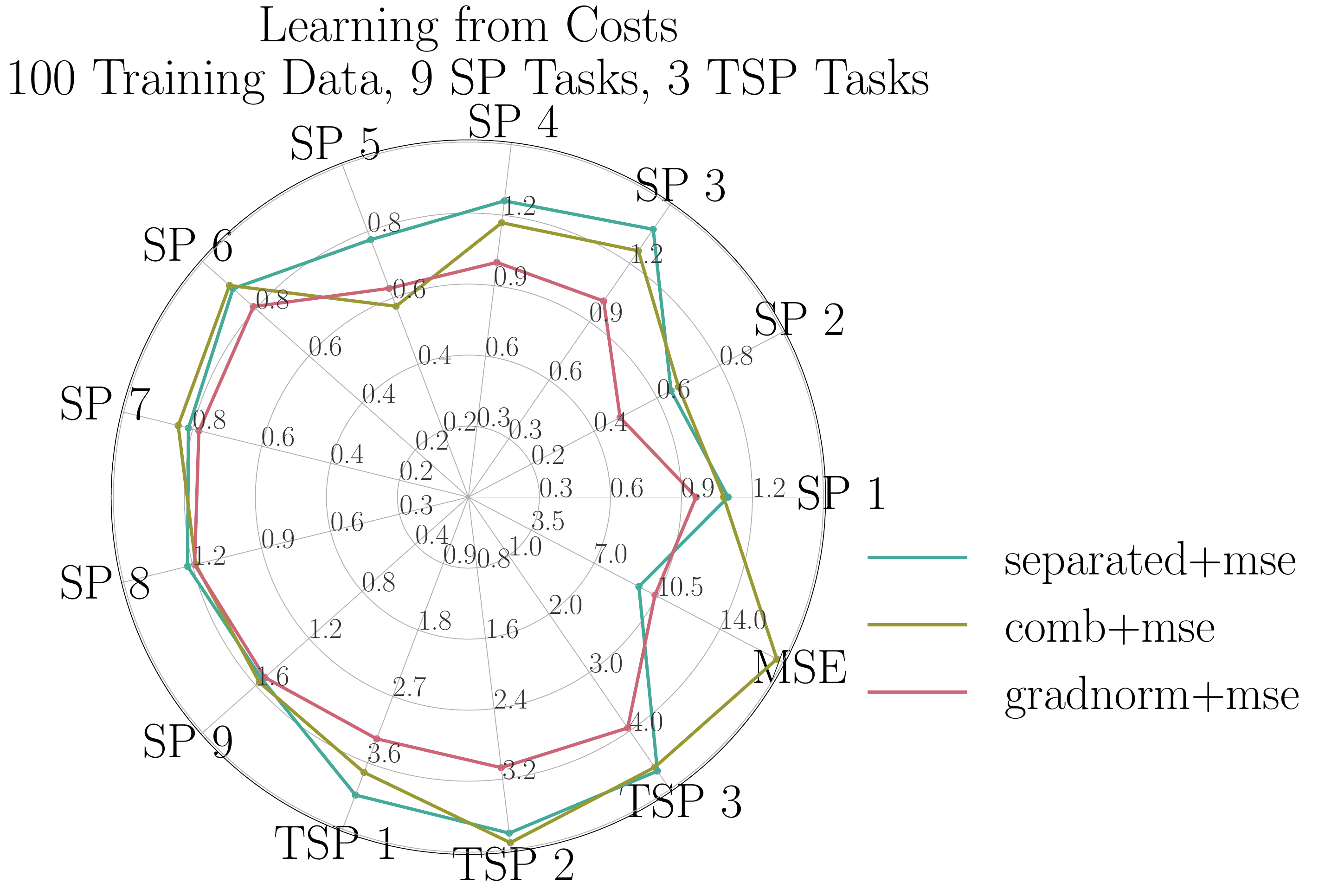}} \\
    \subfloat{\includegraphics[width=0.46\textwidth]{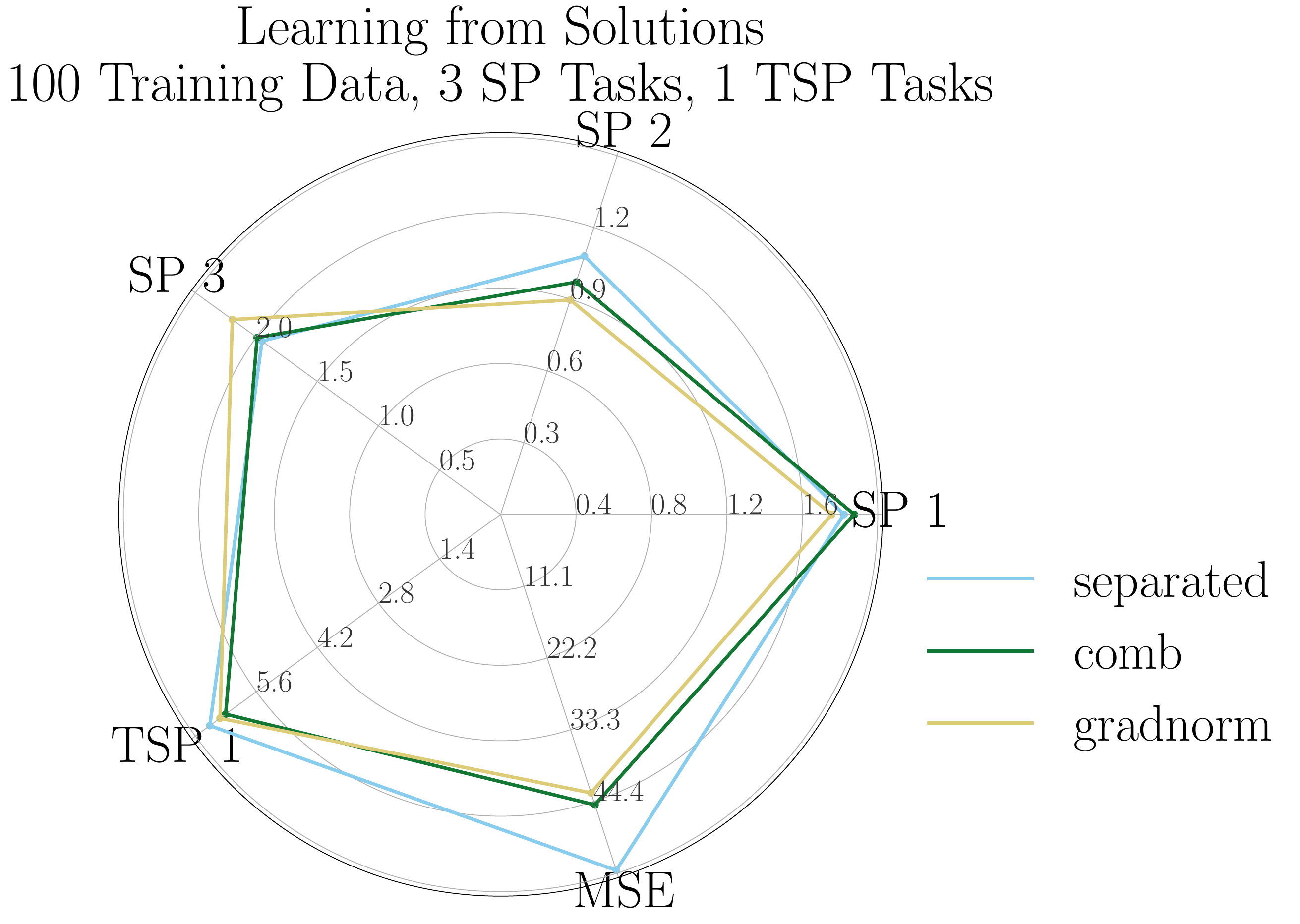}} \hspace{0.3cm}
    \subfloat{\includegraphics[width=0.46\textwidth]{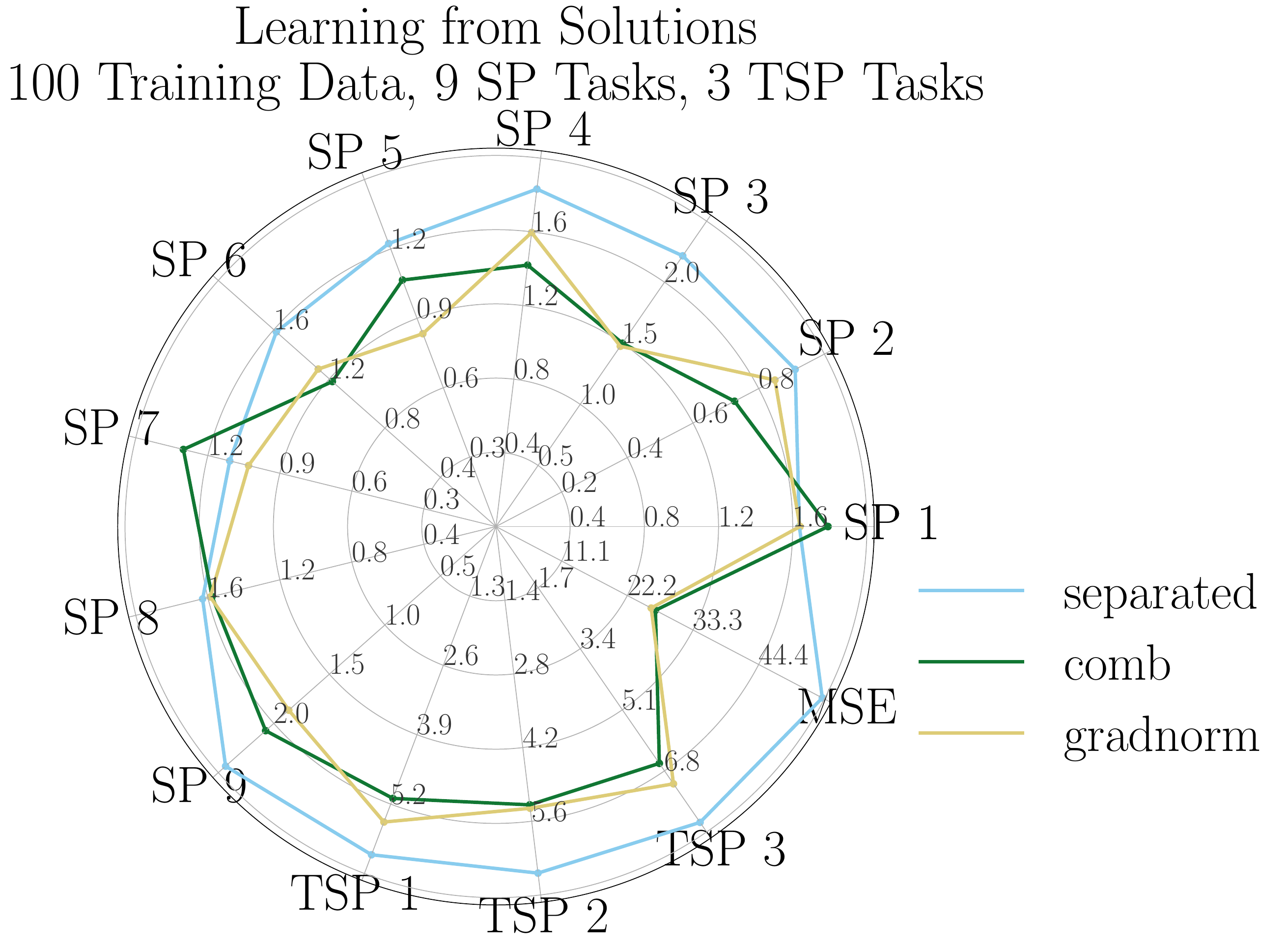}} \\
    \caption{Performance Radar Plot for Graph Routing: Average performance for different tasks on the test set, trained with SPO+ (Top) and PFYL (Bottom), and $100$ training data points, including regrets and cost MSE, lower is better. SP $i$ is the regret for shortest path task $i$,  TSP $i$ is the regret for traveling salesperson task $i$, MSE is the mean squared error of cost coefficients.}
    \label{fig:radar}
\end{figure}

Experimental results on multiple routing tasks on a graph, as shown in Figure \ref{fig:radar}, demonstrate that multi-task end-to-end predict-then-optimize, especially with GradNorm, has a performance advantage over single-task approaches. Figure \ref{fig:radar} shows the results of learning from costs with SPO+ (top), and the results of learning from solutions with PFYL (bottom). In these ``radar plots", the per-task regret on unseen test data is shown along each dimension (lower is better). It can be seen that the innermost method in these figures (best seen in color) is the red one, ``gradnorm+mse". 

\begin{figure}[htbp]
    \centering
    \subfloat{\includegraphics[width=0.50\textwidth]{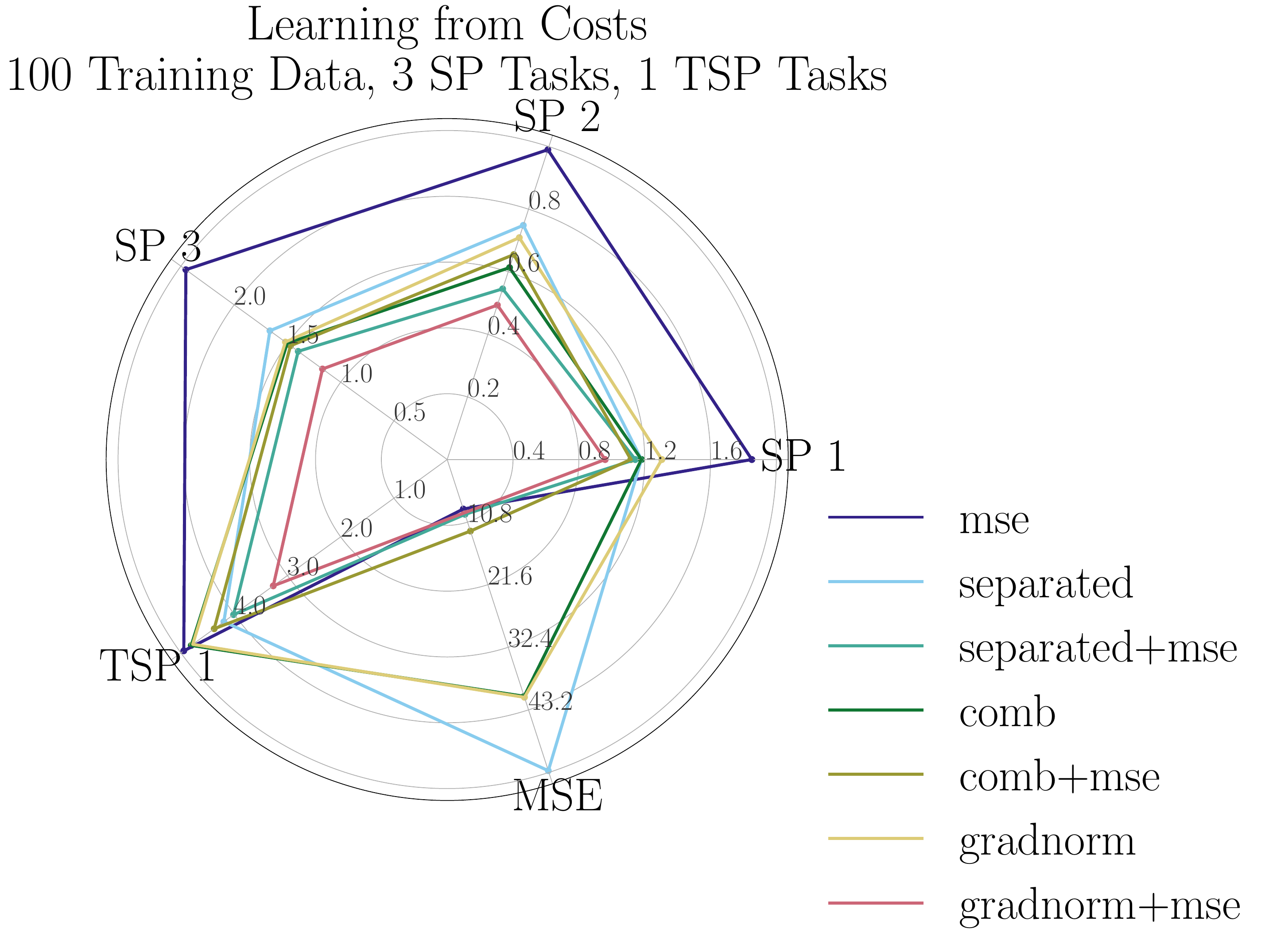}} \hspace{-0.2cm}
    \subfloat{\includegraphics[width=0.50\textwidth]{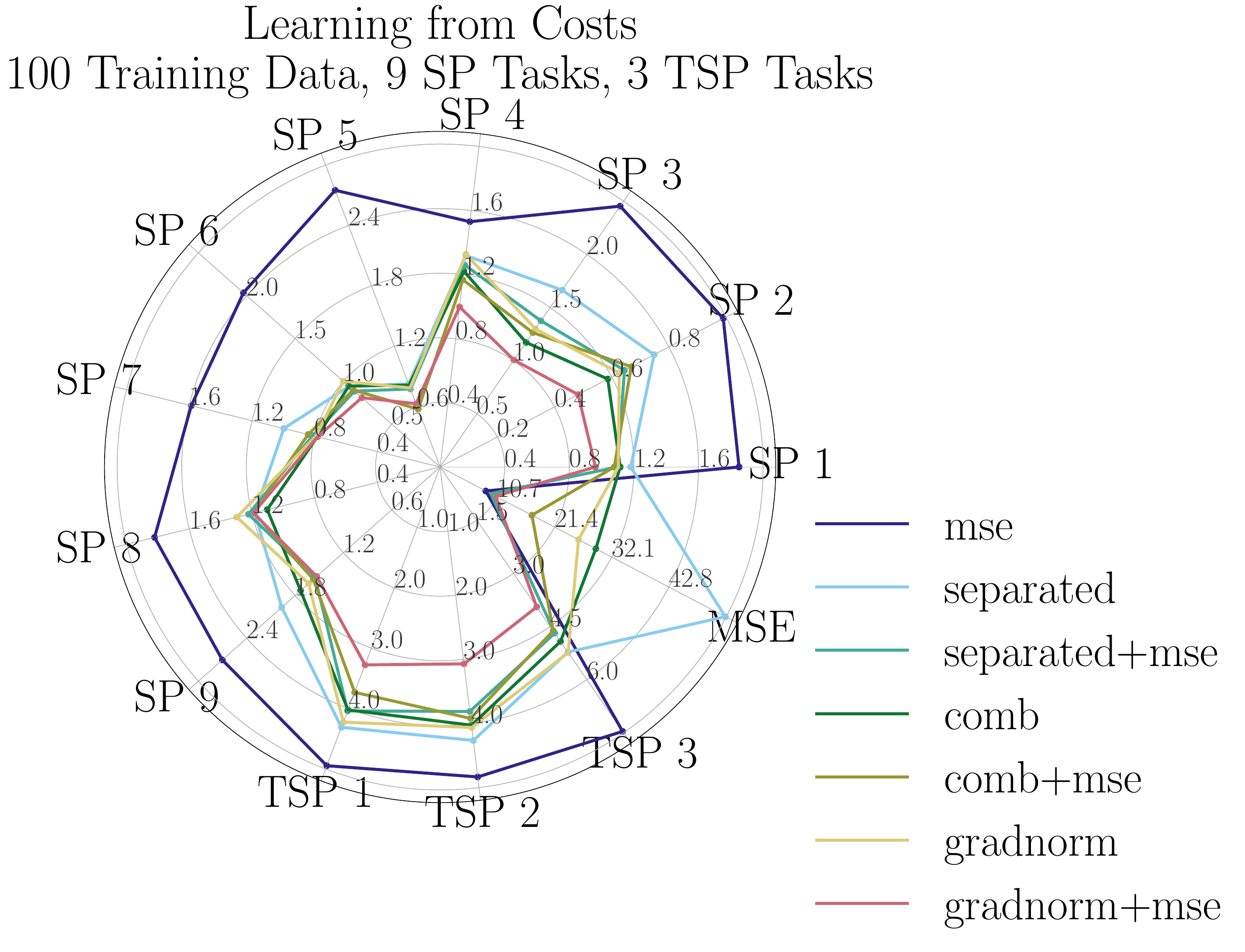}} \\
    \caption{Learning from Costs for Graph Routing: Average performance for different tasks on the test set, trained with SPO+ and $100$ training data points for more methods, lower is better. SP $i$ is the regret for shortest path task $i$,  TSP $i$ is the regret for traveling salesperson task $i$, MSE is the mean squared error of cost coefficients.} 
    \label{fig:spo_radarf}
\end{figure}

More experiments are shown in Fig.~\ref{fig:spo_radarf}. The investigation includes the two-stage method, single-task, and multi-task with and without cost MSE as regularization; these two plots include a superset of the methods in Fig.~\ref{fig:radar}. Despite achieving a lower MSE, the two-stage approaches exhibit significantly worse regret than end-to-end learning. Additionally, adding an MSE regularizer on the cost coefficients consistently improves end-to-end learning approaches. Thus, we always include the additional cost of MSE regularizer when learning from costs. However, since labeled costs are absent when learning directly from solutions, PFYL cannot add cost MSE (``+mse") as regularization and cannot be compared with the two-stage method (which requires cost labels.)

\subsection{Efficiency Benefit of Multi-Task Learning}
\label{subsec:effi}

\begin{figure}[htbp]
    \centering
    \subfloat{\includegraphics[width=0.50\textwidth]{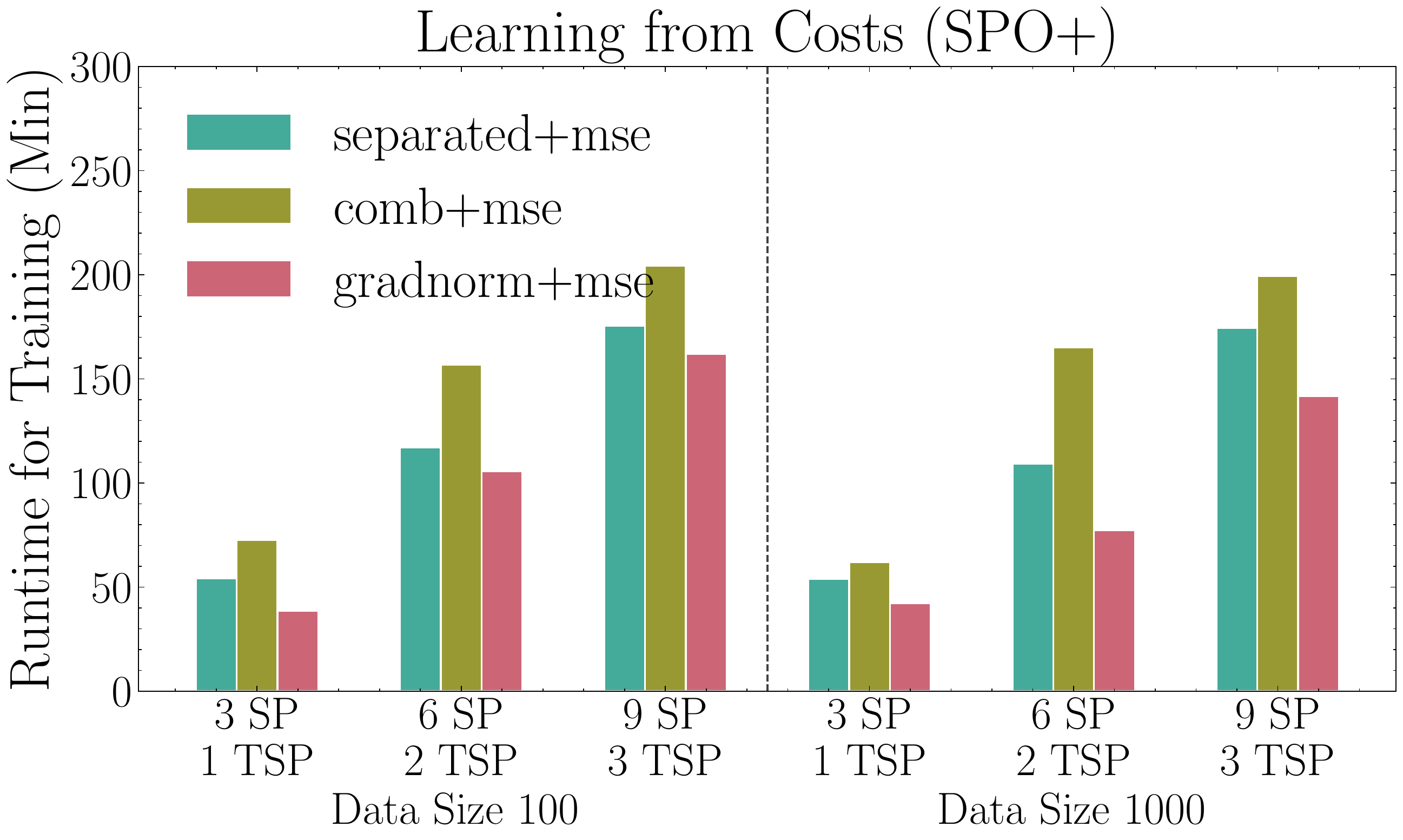}} 
    \subfloat{\includegraphics[width=0.50\textwidth]{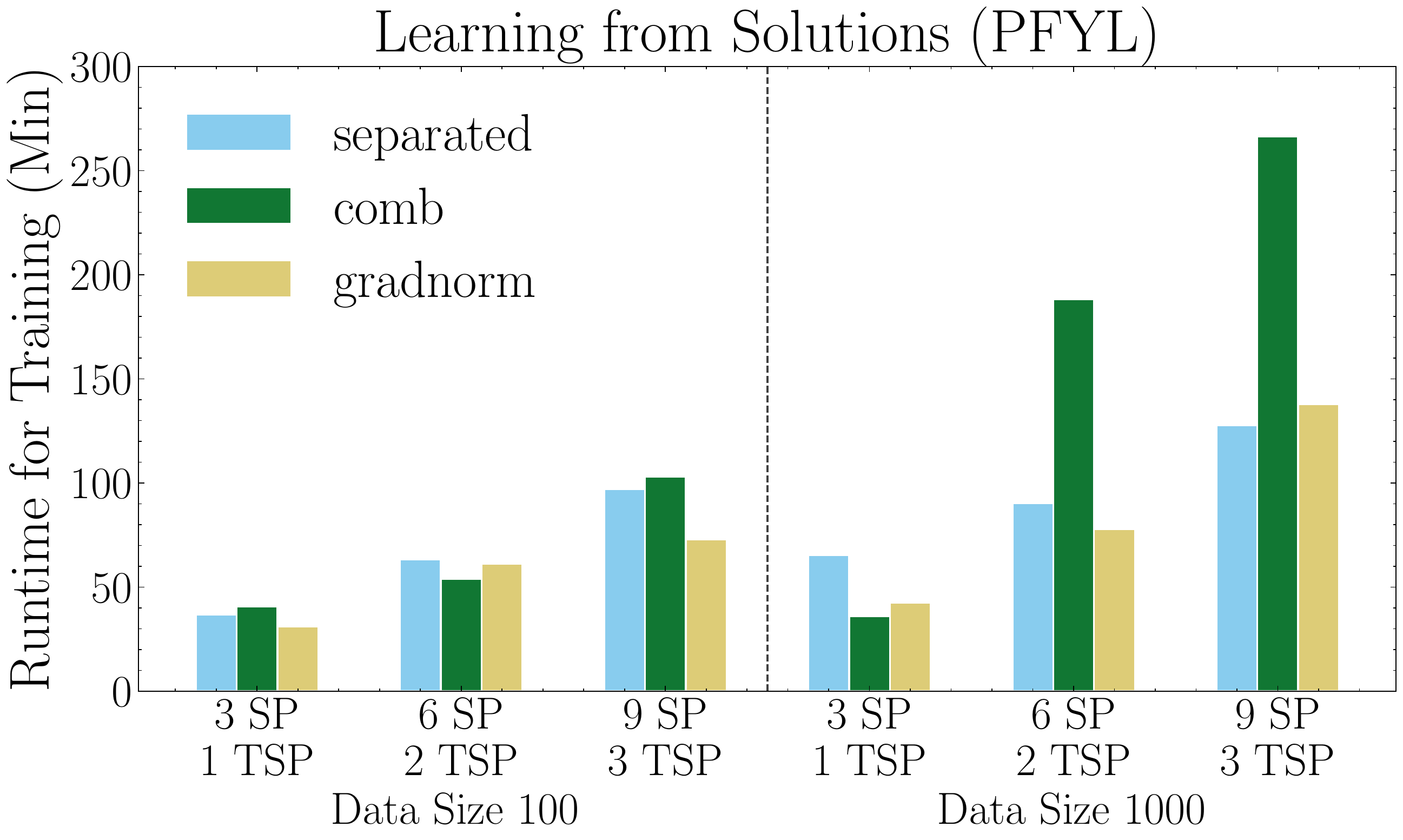}} \\
    \caption{Training time for Graph Routing: The elapsed time of training to convergence at different settings.}
    \label{fig:time}
\end{figure}

Figure \ref{fig:time} shows the training time for SPO+ and PFYL models when using early stopping when five consecutive epochs exhibit non-improving loss values on held-out validation data, a standard trick in neural network training. For ``separated" and ``separated+mse", the training time is the sum of each individual model. We can see that the use of GradNorm to adjust the weights dynamically allows for efficient model training as faster convergence is achieved. The ``separated+mse" baseline typically requires more time to converge, but also converges to worse models in terms of regret, as seen in the previous paragraph. Furthermore, ``comb" and ``comb+mse" usually require more time to converge, which highlights the importance of an effective weighting scheme for multi-tasks approaches.

\subsection{Learning under Data Scarcity}
\label{subsec:lessdata}

\begin{figure}[htbp]
    \hspace{-0.5cm}
    \centering
    \subfloat{\includegraphics[width=0.50\textwidth]{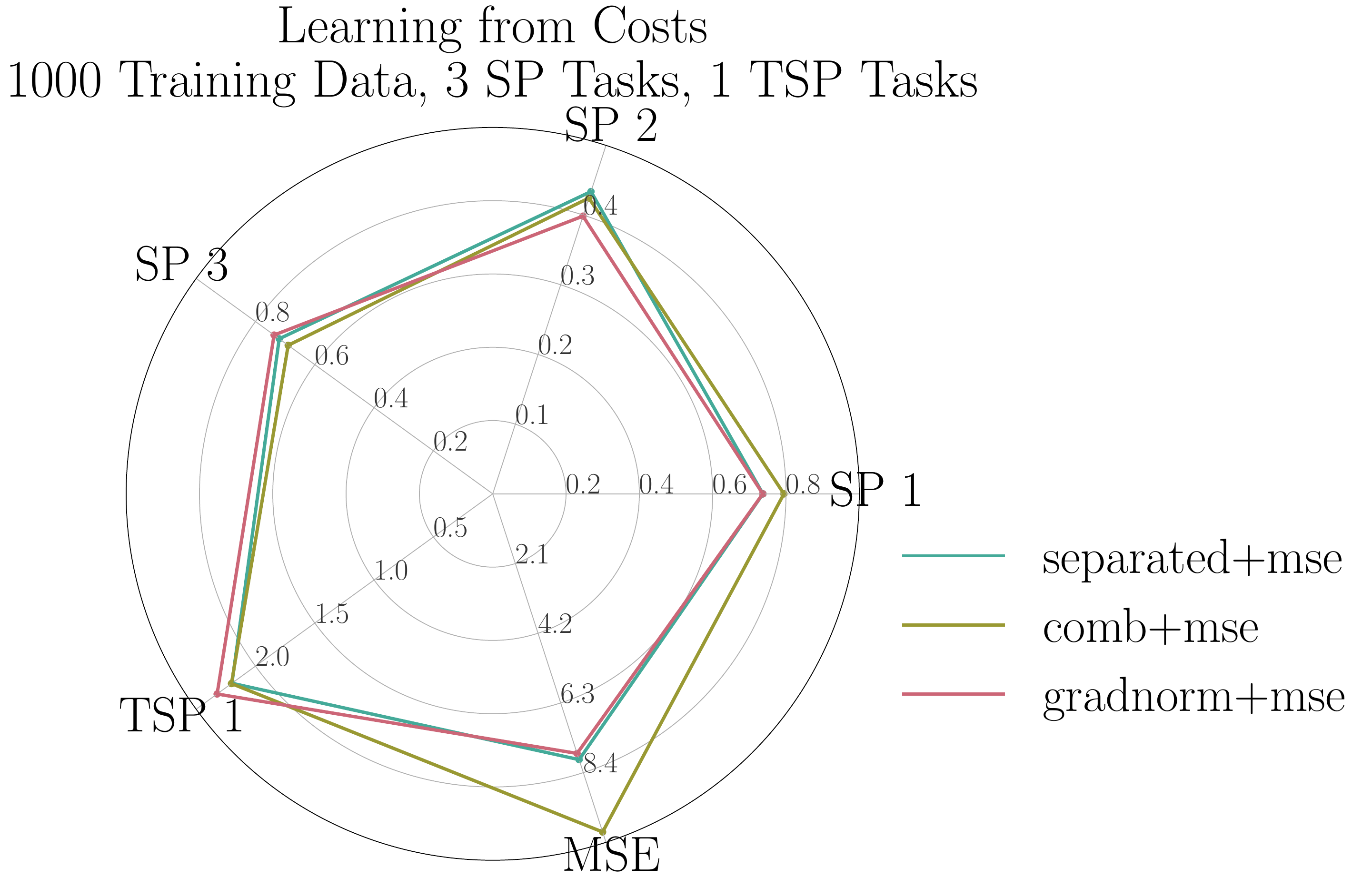}} 
    \subfloat{\includegraphics[width=0.46\textwidth]{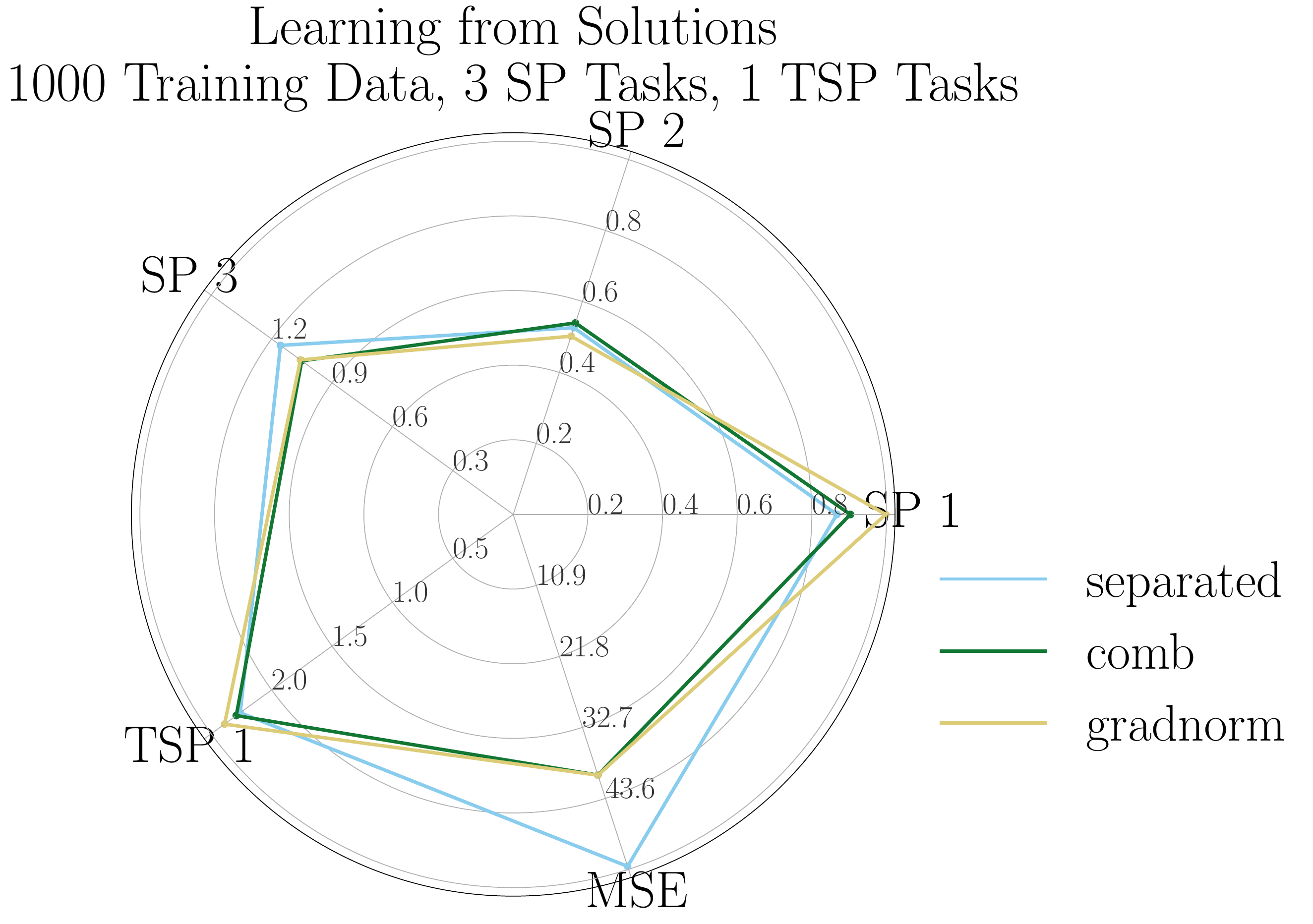}} \\
    \caption{More Training Data for Graph Routing: Average performance for different graph routing tasks on the test set for SPO+ (left) and PFYL (right), trained with PFYL and $1000$ training data points, including regrets and cost MSE, lower is better. SP $i$ is the regret for shortest path task $i$,  TSP $i$ is the regret for traveling salesperson task $i$, MSE is the mean squared error of cost coefficients.}
    \label{fig:moredata}
\end{figure}

\begin{figure}[htbp]
    \centering
    \subfloat{\includegraphics[width=0.33\textwidth]{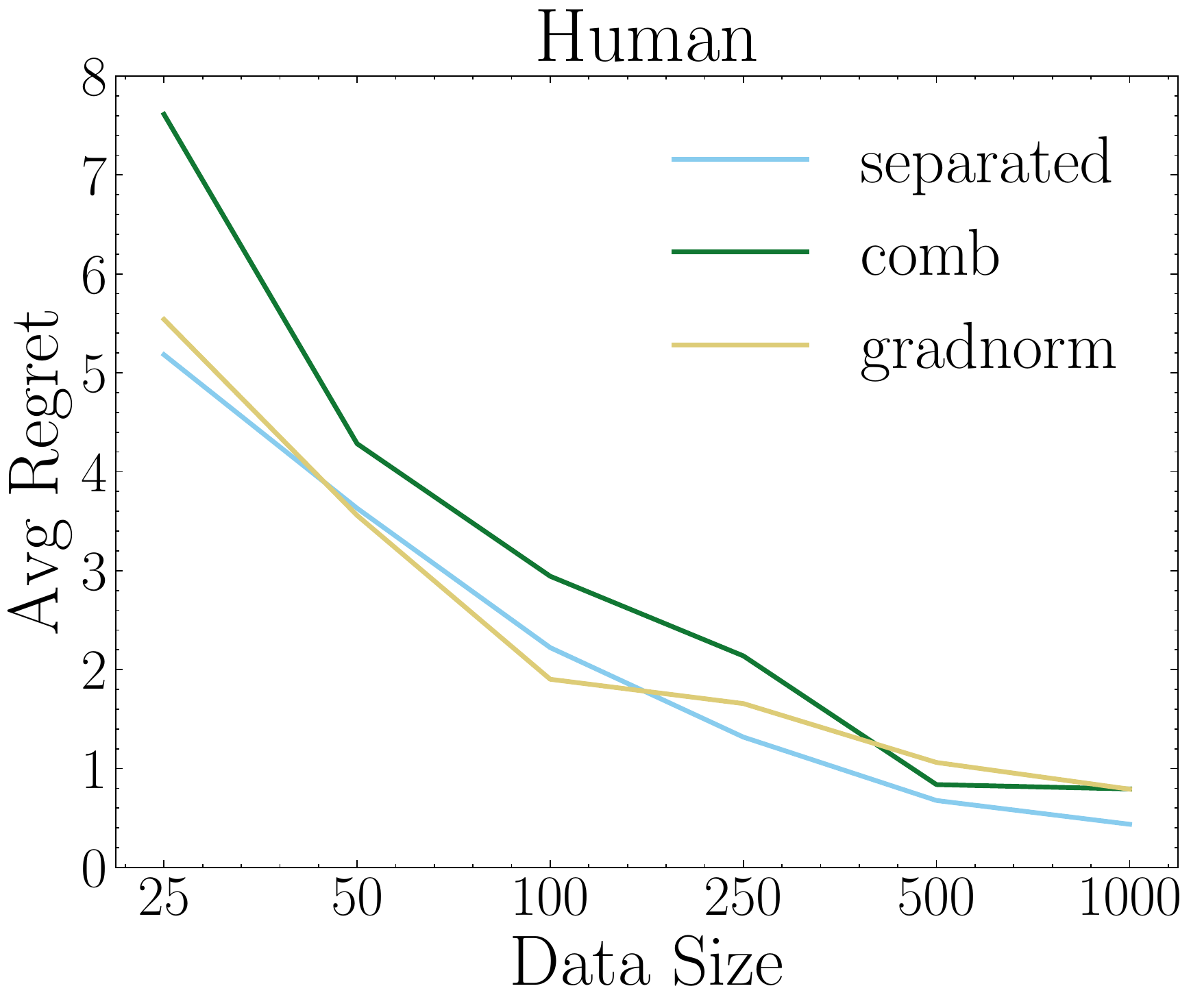}} 
    \subfloat{\includegraphics[width=0.33\textwidth]{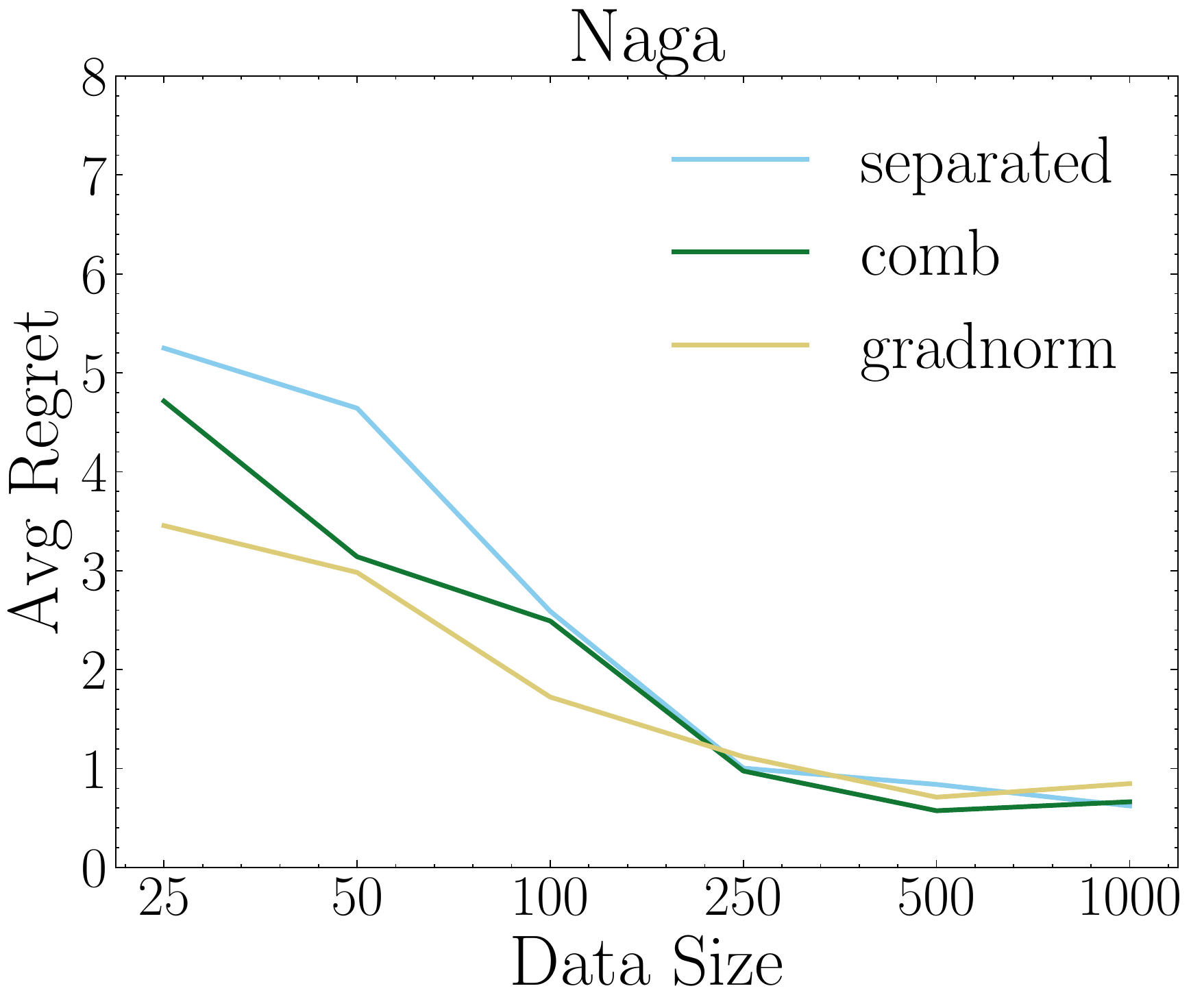}}
    \subfloat{\includegraphics[width=0.33\textwidth]{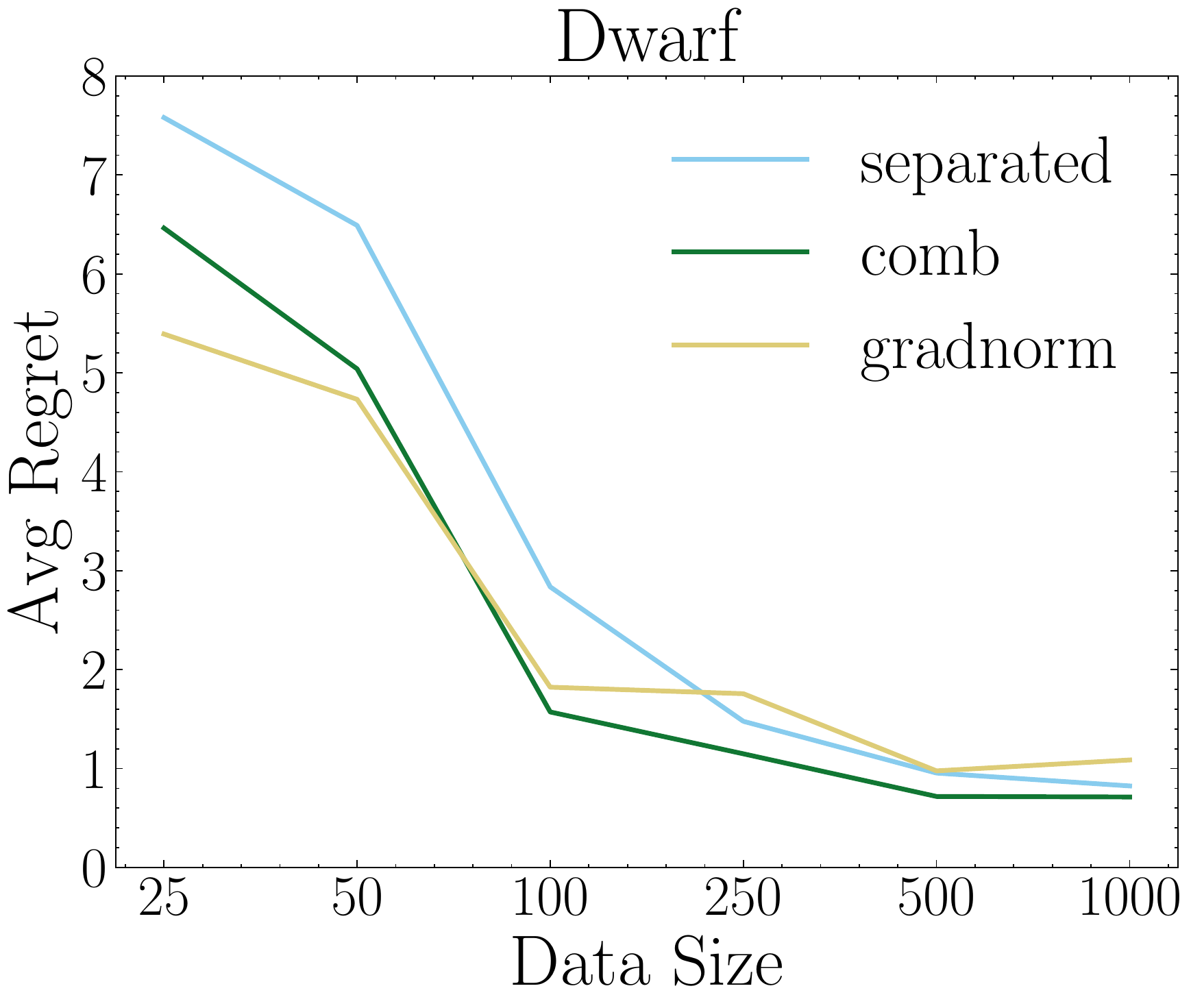}}\\
    \caption{Performance on Warcraft Shortest Path: Average regrets of different strategies, trained with PFYL, decreases as the amount of training data increases; lower is better.}
    \label{fig:data}
\end{figure}

In this section, we claim that the multi-task end-to-end predict-then-optimize framework is particularly effective in the small-data regime. Compared to Figures \ref{fig:radar}, we find that multi-task learning for graph routing loses its advantage with more training data (Figure \ref{fig:moredata}). In Warcraft shortest path problem, Figure \ref{fig:data} shows that the performance of the separated single-task model gradually improves and may even surpass multi-task learning as the amount of training data increases. These figures show that multi-task end-to-end predict-then-optimize can effectively leverage information from related datasets when the size of the individual dataset is limited. Therefore, multi-task learning is a reliable option under data scarcity.

\subsection{Learning under Task Redundancy}
\label{subsec:moretasks}

\begin{figure}[htbp]
    \centering
    \subfloat{\includegraphics[width=0.42\textwidth]{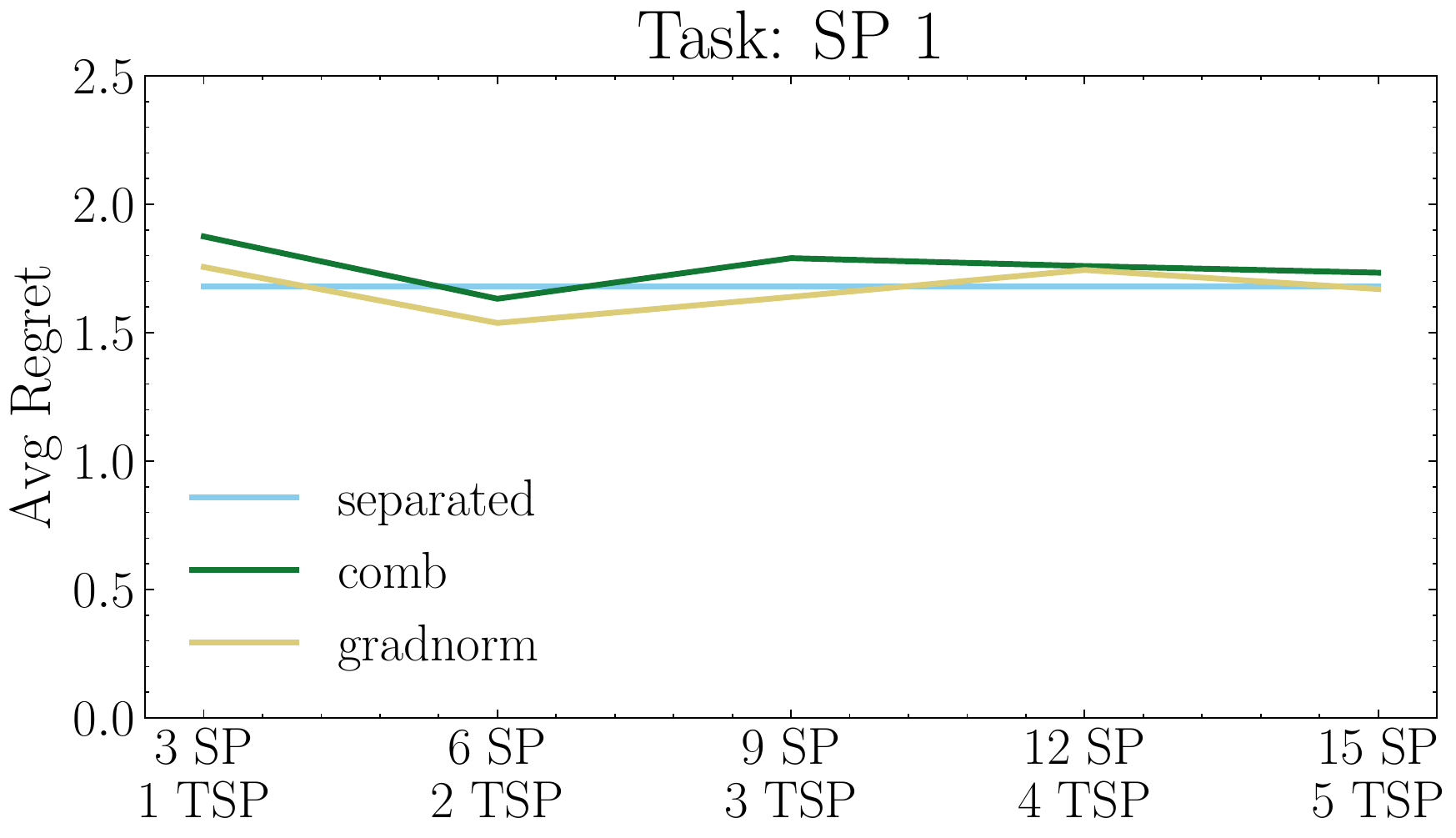}} 
    \subfloat{\includegraphics[width=0.42\textwidth]{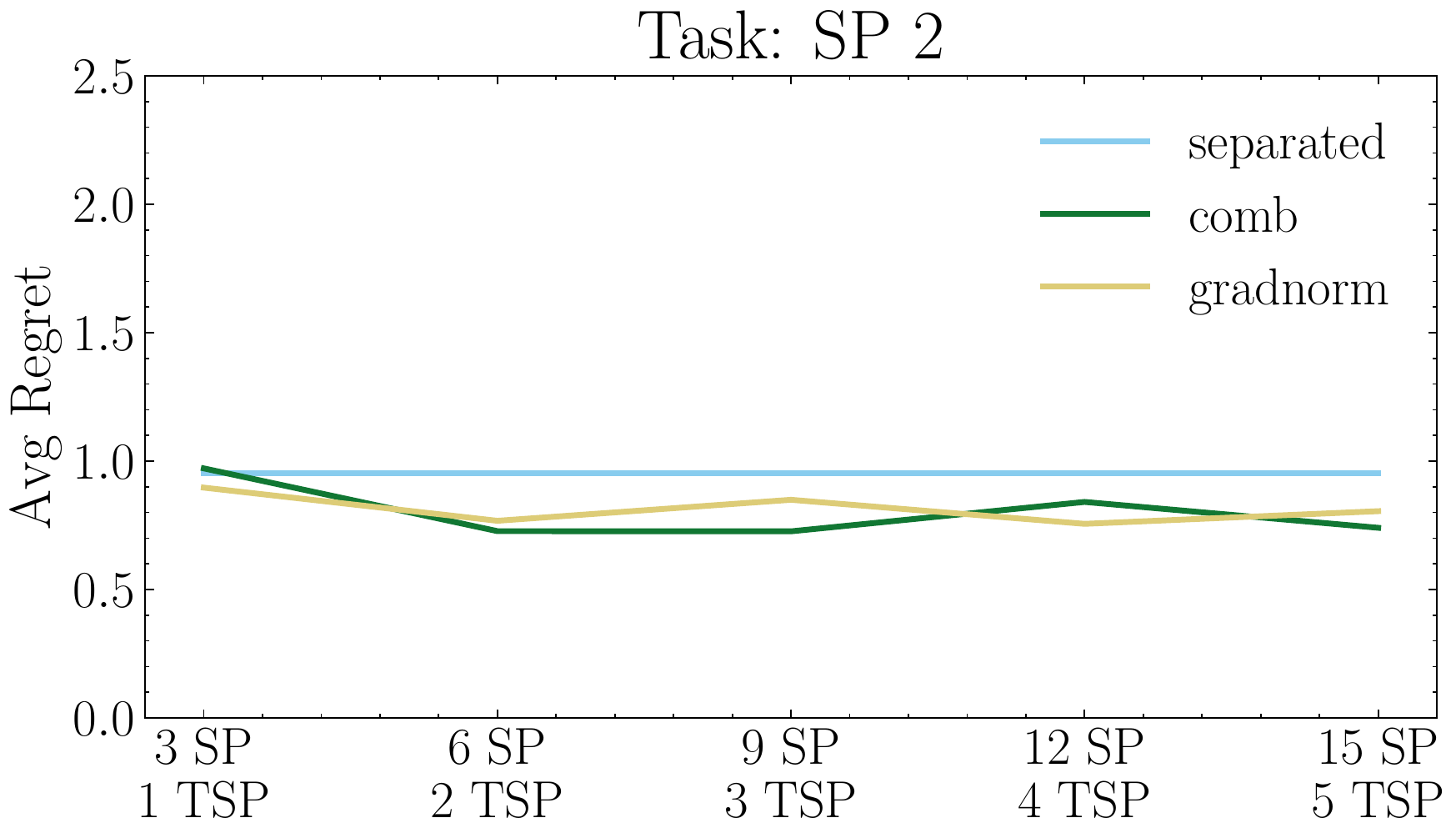}} \\
    \subfloat{\includegraphics[width=0.42\textwidth]{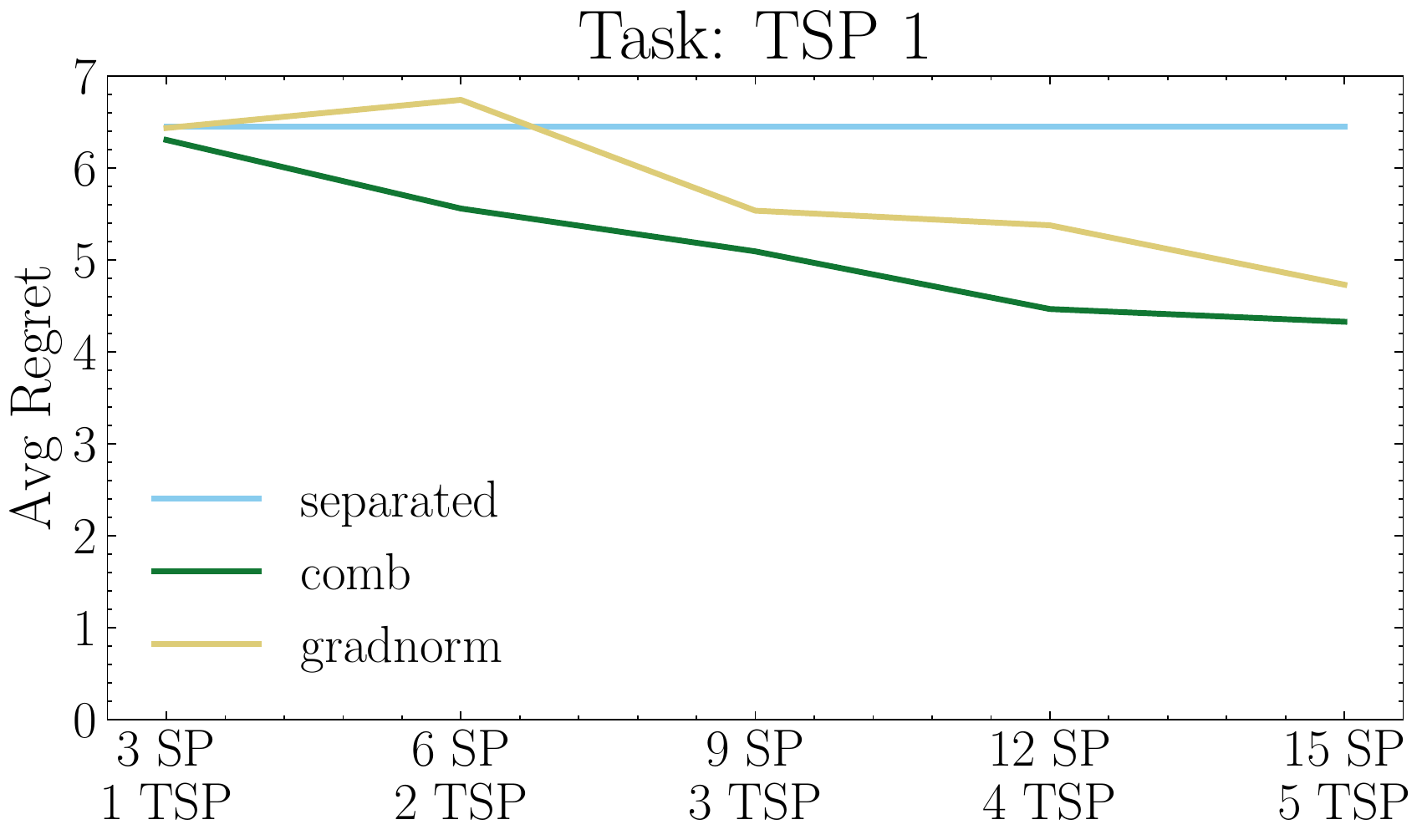}}
    \subfloat{\includegraphics[width=0.42\textwidth]{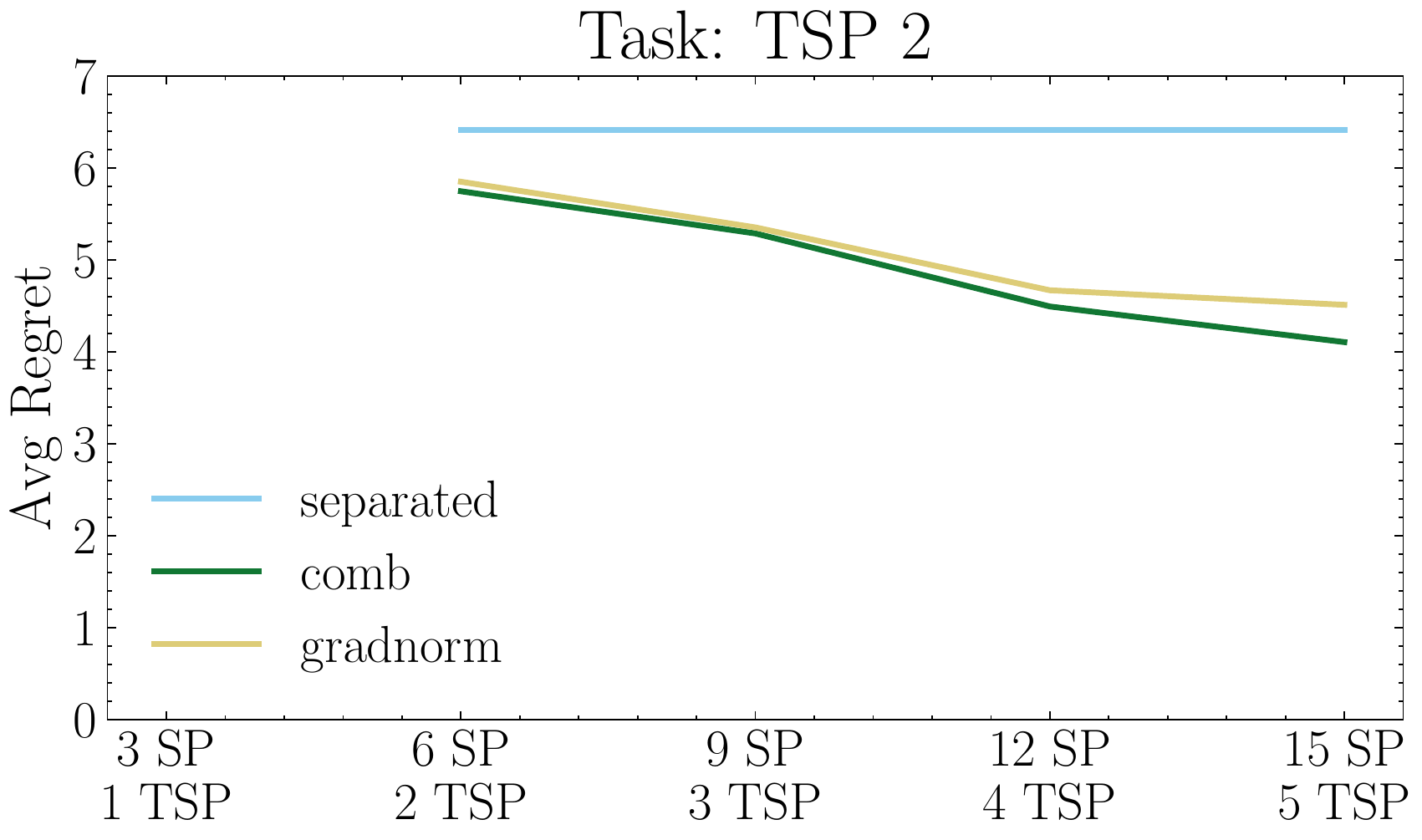}}
    \caption{Performance on Graph Routing: Average regrets of a different strategy, trained with PFYL and $100$ training data, decreases as the amount of tasks increases, lower is better.}
    \label{fig:task}
\end{figure}

Figure \ref{fig:task} indicates that increasing the number of related tasks improves the model performance, especially for complicated tasks such as TSP. This performance improvement can be attributed to the positive interaction between the losses of the related tasks. This finding suggests the potential for using auxiliary tasks to enhance model performance.
\section{Conclusion} \label{sec:conc}

We extend the end-to-end predict-then-optimize framework to multi-task learning, which jointly minimizes decision error for related optimization tasks. Our results demonstrate the benefits of this approach, including an improved performance with less training data points and the ability to handle multiple tasks simultaneously. Future work in this area could include the application of this method to real-world problems, as well as further exploration of techniques for multi-task learning, such as current and novel multi-task neural network architectures and gradient calibration methods.
\\
\\
\\
\\
%
%
%

\let\clearpage\relax
\bibliographystyle{splncs04nat}
\bibliography{ref}

\begin{thebibliography}{34}
\providecommand{\natexlab}[1]{#1}
\providecommand{\url}[1]{\texttt{#1}}
\providecommand{\urlprefix}{URL }
\expandafter\ifx\csname urlstyle\endcsname\relax
  \providecommand{\doi}[1]{doi:\discretionary{}{}{}#1}\else
  \providecommand{\doi}{doi:\discretionary{}{}{}\begingroup
  \urlstyle{rm}\Url}\fi

\bibitem[{Agrawal et~al.(2019)Agrawal, Amos, Barratt, Boyd, Diamond, and
  Kolter}]{agrawal2019differentiable}
Agrawal, A., Amos, B., Barratt, S., Boyd, S., Diamond, S., Kolter, J.Z.:
  Differentiable convex optimization layers. In: Wallach, H., Larochelle, H.,
  Beygelzimer, A., d\textquotesingle Alch\'{e}-Buc, F., Fox, E., Garnett, R.
  (eds.) Advances in Neural Information Processing Systems, vol.~32, Curran
  Associates, Inc. (2019)

\bibitem[{Ahuja and Orlin(2001)}]{ahuja2001inverse}
Ahuja, R.K., Orlin, J.B.: Inverse optimization. Operations Research
  \textbf{49}(5), 771--783 (2001)

\bibitem[{Amos and Kolter(2017)}]{amos2017optnet}
Amos, B., Kolter, J.Z.: Optnet: Differentiable optimization as a layer in
  neural networks. In: International Conference on Machine Learning, pp.
  136--145, PMLR (2017)

\bibitem[{Bengio(1997)}]{bengio1997using}
Bengio, Y.: Using a financial training criterion rather than a prediction
  criterion. International Journal of Neural Systems \textbf{8}(04), 433--443
  (1997)

\bibitem[{Berthet et~al.(2020)Berthet, Blondel, Teboul, Cuturi, Vert, and
  Bach}]{berthet2020learning}
Berthet, Q., Blondel, M., Teboul, O., Cuturi, M., Vert, J.P., Bach, F.:
  Learning with differentiable perturbed optimizers. arXiv preprint
  arXiv:2002.08676  (2020)

\bibitem[{Caruana(1997)}]{caruana1997multitask}
Caruana, R.: Multitask learning. Machine learning \textbf{28}(1), 41--75 (1997)

\bibitem[{Chen et~al.(2018)Chen, Badrinarayanan, Lee, and
  Rabinovich}]{chen2018gradnorm}
Chen, Z., Badrinarayanan, V., Lee, C.Y., Rabinovich, A.: Gradnorm: Gradient
  normalization for adaptive loss balancing in deep multitask networks. In:
  International conference on machine learning, pp. 794--803, PMLR (2018)

\bibitem[{Dalle et~al.(2022)Dalle, Baty, Bouvier, and
  Parmentier}]{dalle2022learning}
Dalle, G., Baty, L., Bouvier, L., Parmentier, A.: Learning with combinatorial
  optimization layers: a probabilistic approach. arXiv preprint
  arXiv:2207.13513  (2022)

\bibitem[{Dantzig et~al.(1954)Dantzig, Fulkerson, and
  Johnson}]{dantzig1954solution}
Dantzig, G., Fulkerson, R., Johnson, S.: Solution of a large-scale
  traveling-salesman problem. Journal of the operations research society of
  America \textbf{2}(4), 393--410 (1954)

\bibitem[{Donti et~al.(2017)Donti, Amos, and Kolter}]{donti2017task}
Donti, P.L., Amos, B., Kolter, J.Z.: Task-based end-to-end model learning in
  stochastic optimization. Advances in Neural Information Processing Systems
  (2017)

\bibitem[{Duong et~al.(2015)Duong, Cohn, Bird, and Cook}]{duong2015low}
Duong, L., Cohn, T., Bird, S., Cook, P.: Low resource dependency parsing:
  Cross-lingual parameter sharing in a neural network parser. In: Proceedings
  of the 53rd annual meeting of the Association for Computational Linguistics
  and the 7th international joint conference on natural language processing
  (volume 2: short papers), pp. 845--850 (2015)

\bibitem[{Elmachtoub and Grigas(2021)}]{elmachtoub2021smart}
Elmachtoub, A.N., Grigas, P.: Smart “predict, then optimize”. Management
  Science  (2021)

\bibitem[{Estes and Richard(2023)}]{estes2023smart}
Estes, A.S., Richard, J.P.P.: Smart predict-then-optimize for two-stage linear
  programs with side information. INFORMS Journal on Optimization  (2023)

\bibitem[{Ferber et~al.(2020)Ferber, Wilder, Dilkina, and
  Tambe}]{ferber2020mipaal}
Ferber, A., Wilder, B., Dilkina, B., Tambe, M.: Mipaal: Mixed integer program
  as a layer. In: Proceedings of the AAAI Conference on Artificial
  Intelligence, vol.~34, pp. 1504--1511 (2020)

\bibitem[{Ford et~al.(2015)Ford, Nguyen, Tambe, Sintov, and
  Delle~Fave}]{ford2015beware}
Ford, B., Nguyen, T., Tambe, M., Sintov, N., Delle~Fave, F.: Beware the
  soothsayer: From attack prediction accuracy to predictive reliability in
  security games. In: International Conference on Decision and Game Theory for
  Security, pp. 35--56, Springer (2015)

\bibitem[{{Gurobi Optimization, LLC}(2021)}]{gurobi}
{Gurobi Optimization, LLC}: {Gurobi Optimizer Reference Manual} (2021),
  \urlprefix\url{https://www.gurobi.com}

\bibitem[{Hu et~al.(2023)Hu, Lee, and Lee}]{hu2023predict}
Hu, X., Lee, J.C., Lee, J.H.: Predict+ optimize for packing and covering lps
  with unknown parameters in constraints. In: Proceedings of the AAAI
  Conference on Artificial Intelligence (2023)

\bibitem[{Kendall et~al.(2018)Kendall, Gal, and Cipolla}]{kendall2018multi}
Kendall, A., Gal, Y., Cipolla, R.: Multi-task learning using uncertainty to
  weigh losses for scene geometry and semantics. In: Proceedings of the IEEE
  conference on computer vision and pattern recognition, pp. 7482--7491 (2018)

\bibitem[{Liu et~al.(2021)Liu, Liu, Jin, Stone, and Liu}]{liu2021conflict}
Liu, B., Liu, X., Jin, X., Stone, P., Liu, Q.: Conflict-averse gradient descent
  for multi-task learning. Advances in Neural Information Processing Systems
  \textbf{34}, 18878--18890 (2021)

\bibitem[{Liu et~al.(2019)Liu, Johns, and Davison}]{liu2019end}
Liu, S., Johns, E., Davison, A.J.: End-to-end multi-task learning with
  attention. In: Proceedings of the IEEE/CVF conference on computer vision and
  pattern recognition, pp. 1871--1880 (2019)

\bibitem[{Ma et~al.(2018)Ma, Zhao, Yi, Chen, Hong, and Chi}]{ma2018modeling}
Ma, J., Zhao, Z., Yi, X., Chen, J., Hong, L., Chi, E.H.: Modeling task
  relationships in multi-task learning with multi-gate mixture-of-experts. In:
  Proceedings of the 24th ACM SIGKDD international conference on knowledge
  discovery \& data mining, pp. 1930--1939 (2018)

\bibitem[{Mandi and Guns(2020)}]{mandi2020interior}
Mandi, J., Guns, T.: Interior point solving for lp-based
  prediction+optimisation. In: Larochelle, H., Ranzato, M., Hadsell, R.,
  Balcan, M.F., Lin, H. (eds.) Advances in Neural Information Processing
  Systems, vol.~33, pp. 7272--7282, Curran Associates, Inc. (2020)

\bibitem[{Mandi et~al.(2020)Mandi, Stuckey, Guns et~al.}]{mandi2020smart}
Mandi, J., Stuckey, P.J., Guns, T., et~al.: Smart predict-and-optimize for hard
  combinatorial optimization problems. In: Proceedings of the AAAI Conference
  on Artificial Intelligence, vol.~34, pp. 1603--1610 (2020)

\bibitem[{Mulamba et~al.(2020)Mulamba, Mandi, Diligenti, Lombardi, Bucarey, and
  Guns}]{mulamba2020contrastive}
Mulamba, M., Mandi, J., Diligenti, M., Lombardi, M., Bucarey, V., Guns, T.:
  Contrastive losses and solution caching for predict-and-optimize. arXiv
  preprint arXiv:2011.05354  (2020)

\bibitem[{Paszke et~al.(2019)Paszke, Gross, Massa, Lerer, Bradbury, Chanan,
  Killeen, Lin, Gimelshein, Antiga et~al.}]{paszke2019pytorch}
Paszke, A., Gross, S., Massa, F., Lerer, A., Bradbury, J., Chanan, G., Killeen,
  T., Lin, Z., Gimelshein, N., Antiga, L., et~al.: Pytorch: An imperative
  style, high-performance deep learning library. Advances in neural information
  processing systems \textbf{32} (2019)

\bibitem[{Shah et~al.(2022)Shah, Wang, Wilder, Perrault, and
  Tambe}]{shahdecision}
Shah, S., Wang, K., Wilder, B., Perrault, A., Tambe, M.: Decision-focused
  learning without decision-making: Learning locally optimized decision losses.
  In: Advances in Neural Information Processing Systems (2022)

\bibitem[{Shazeer et~al.(2017)Shazeer, Mirhoseini, Maziarz, Davis, Le, Hinton,
  and Dean}]{shazeer2017outrageously}
Shazeer, N., Mirhoseini, A., Maziarz, K., Davis, A., Le, Q., Hinton, G., Dean,
  J.: Outrageously large neural networks: The sparsely-gated mixture-of-experts
  layer. arXiv preprint arXiv:1701.06538  (2017)

\bibitem[{Tang and Khalil(sion)}]{tang2022pyepo}
Tang, B., Khalil, E.B.: Pyepo: A pytorch-based end-to-end predict-then-optimize
  library for linear and integer programming. Mathematical Programming
  Computation  (2022 in submission)

\bibitem[{Tang et~al.(2020)Tang, Liu, Zhao, and Gong}]{tang2020progressive}
Tang, H., Liu, J., Zhao, M., Gong, X.: Progressive layered extraction (ple): A
  novel multi-task learning (mtl) model for personalized recommendations. In:
  Fourteenth ACM Conference on Recommender Systems, pp. 269--278 (2020)

\bibitem[{Vlastelica et~al.(2019)Vlastelica, Paulus, Musil, Martius, and
  Rol{\'\i}nek}]{poganvcic2019differentiation}
Vlastelica, M., Paulus, A., Musil, V., Martius, G., Rol{\'\i}nek, M.:
  Differentiation of blackbox combinatorial solvers. arXiv preprint
  arXiv:1912.02175  (2019)

\bibitem[{Wang et~al.(2020)Wang, Tsvetkov, Firat, and Cao}]{wang2020gradient}
Wang, Z., Tsvetkov, Y., Firat, O., Cao, Y.: Gradient vaccine: Investigating and
  improving multi-task optimization in massively multilingual models. arXiv
  preprint arXiv:2010.05874  (2020)

\bibitem[{Wilder et~al.(2019)Wilder, Dilkina, and Tambe}]{wilder2019melding}
Wilder, B., Dilkina, B., Tambe, M.: Melding the data-decisions pipeline:
  Decision-focused learning for combinatorial optimization. In: Proceedings of
  the AAAI Conference on Artificial Intelligence, vol.~33, pp. 1658--1665
  (2019)

\bibitem[{Winkenbach et~al.(2021)Winkenbach, Parks, and
  Noszek}]{winkenbach2021technical}
Winkenbach, M., Parks, S., Noszek, J.: Technical proceedings of the amazon last
  mile routing research challenge  (2021)

\bibitem[{Yu et~al.(2020)Yu, Kumar, Gupta, Levine, Hausman, and
  Finn}]{yu2020gradient}
Yu, T., Kumar, S., Gupta, A., Levine, S., Hausman, K., Finn, C.: Gradient
  surgery for multi-task learning. Advances in Neural Information Processing
  Systems \textbf{33}, 5824--5836 (2020)

\end{thebibliography}



\end{document}